\documentclass[11pt]{amsart}
\usepackage{epsfig}
\usepackage{amsthm}
\usepackage{amssymb}
\usepackage{amsmath}
\usepackage{epic}
\usepackage{eepic}

\newcommand     {\comment}[1]   {}
\newcommand{\mute}[2] {}
\newcommand     {\printname}[1] {}
\newcommand{\labell}[1] {\label{#1}\printname{#1}}

\newtheorem {Theorem}   {Theorem}
\newtheorem {Lemma}[equation]    {Lemma}
\newtheorem {Claim}[equation]    {Claim}
\newtheorem {Proposition}[equation]{Proposition}
\theoremstyle{definition}

\newtheorem {Corollary}[equation]{Corollary}
\theoremstyle{remark}
\newtheorem{Remark}[equation]{Remark}

\def    \valpha      {\alpha}
\def	\vV	{{V}}

\def	\bfC	{{\mathbf C}}
\def	\bfS	{{\mathbf S}}
\def	\bfO	{{\mathbf O}}
\def	\bfM	{{\mathbf M}}

\def    \inv    {^{-1}}

\def    \ssminus        {{\smallsetminus}}

\def    \l<        {\left< }
\def    \r>        {\right> }
\def	\del	{\partial}
\newcommand {\deldel}[1] {\frac{\partial}{\partial #1}}

\def	\nonneg	{{\geq 0}}

\def	\span	{{\operatorname{span}}}

\def	\Td	{{\operatorname{Td}}}

\def	\half	{{\frac{1}{2}}}

\def	\Z	{{\mathbb Z}}
\def	\R	{{\mathbb R}}
\def	\C	{{\mathbb C}}
 \def \L{{\bf L}}
 \def \Ltk{{\bf L}^{2k}}

 \def	\half	{{\frac{1}{2}}}

 \def	\st	{{\text{st}}}

\begin{document}

\title[Euler Maclaurin for a simple integral polytope]
{Euler Maclaurin with remainder for a simple integral polytope}

\author[Y.\ Karshon]{Yael Karshon}
\address{Institute of Mathematics, The Hebrew University of Jerusalem,
Israel, and: Department of Mathematics, The University of Toronto, 
Toronto, Ontario M5S 3G3, Canada}
\email{karshon@math.toronto.edu}

\author[S.\ Sternberg]{Shlomo Sternberg}
\address{Department of Mathematics, Harvard University,
Cambridge, MA 02138, USA}
\email{shlomo@math.harvard.edu}

\author[J.\ Weitsman]{Jonathan Weitsman}
\address{Department of Mathematics, University of California,
Santa Cruz, CA 95064, USA}
\email{weitsman@math.UCSC.EDU}

\thanks{This work was partially supported by
United States -- Israel Binational Science Foundation
grant number 2000352 (to Y.K. and J.W.),  by the Connaught Fund (to Y.K.),
and by National Science Foundation Grant DMS 99/71914 (to J.W.).}

\begin{abstract}
We give an Euler Maclaurin formula with remainder for the sum of the
values of a smooth function on the integral points in a simple integral 
polytope.  We prove this formula by elementary methods.
\end{abstract}
\maketitle
%\tableofcontents

% -------------------------------------------------------------------------
\section{Introduction}
% -------------------------------------------------------------------------
\labell{sec:intro}

The Euler Maclaurin formula computes the sum of the values
of a function $f$ over the integer points in an interval
in terms of the integral of $f$ over variations of that interval.
A version of this classical formula is this: 

For any function $f(x)$ on the real line and any integers $a<b$,
we will consider the weighted sum
\begin{equation} \labell{wted sum}
 \sum_{[a,b]}{'} f := \half f(a) + f(a+1) + \ldots + f(b-1) + \half f(b) . 
\end{equation}
If $f$ is ``nice enough", for instance, a polynomial, then
\begin{equation} \labell{EM exact weighted}
   \sum_{[a,b]}{'} f =
   \left. \L(\deldel{h_1}) \L(\deldel{h_2}) \int_{a-h_1}^{b+h_2} f(x) dx 
   \right|_{h_1 = h_2 = 0},
\end{equation}
where
\begin{equation} \labell{LS}
 \L(S) = \frac{S/2}{\tanh(S/2)} 
       = 1 + \sum_{k=1}^\infty \frac{1}{(2k)!} b_{2k} S^{2k} .
\end{equation}
Because $\int_{a-h_1}^{b+h_2} f(x) dx$ is a polynomial in $h_1$ and $h_2$
if $f$ is a polynomial in $x$,
applying the infinite order differential operator $\L(\deldel{h_i})$
then yields a finite sum, so
the right hand side of \eqref{EM exact weighted} is well defined
when $f$ is a polynomial.

A polytope in $\R^n$ is called \emph{integral}, or a \emph{lattice polytope},
if its vertices
are in the lattice $\Z^n$;
it is called \emph{simple} if exactly $n$ edges emanate from each vertex;
it is called \emph{regular} if, additionally,
the edges emanating from each vertex lie along lines which are generated
by a $\Z$-basis of the lattice $\Z^n$.

Khovanskii and Pukhlikov \cite{KP1,KP2}, following
Khovanskii \cite{Kh1,Kh2},
generalized the classical Euler Maclaurin formula to give
a formula for the sums of the values of polynomial or exponential
functions on the lattice points in higher dimensional 
convex polytopes $\Delta$ which are integral and regular.
This formula was generalized to \emph{simple} integral polytopes 
by Cappell and Shaneson \cite{CS:bulletin,CS:EM,CS:private,S},
and subsequently by Guillemin \cite{Gu} and by Brion-Vergne \cite{BV}.
All of these generalizations involve ``corrections" to the 
Khovanskii-Pukhlikov formula when the simple polytope is not regular.
When applied to the constant function $f \equiv 1$, these Euler Maclaurin
formulas compute the number
of lattice points in $\Delta$ in terms of the volumes
of ``dilations" of $\Delta.$
A small sample of the literature on the problem of counting lattice
points in convex polytopes is given in
\cite{pick,Md,V,KK,Mo,Po,DR,BDR,Ha};
see the survey \cite{BP} and references therein.

These formulas are closely related to the Riemann Roch formula
from algebraic geometry via the correspondence between
polytopes and toric varieties.
Under this correspondence, regular polytopes correspond
to smooth toric varieties, and the Khovanskii-Pukhlikov
formula was motivated in this way, although it was proved combinatorially.
Cappell and Shaneson derived their formula
from their theory of characteristic classes of singular algebraic
varieties and had the key idea of using the operator $\L$ as we do here.
Guillemin obtained his formula from the equivariant
Kawasaki-Riemann-Roch formula and methods coming from
symplectic geometry and the theory of geometric quantization.
Brion and Vergne employed a method 
that is closer to that in the original proof of Khovanskii and
Pukhlikov, using Fourier analysis.

To illustrate the relation to toric varieties, let us sketch
the symplectic-geometric proof for the case of a regular polytope,
following Guillemin \cite{Gu:book}.  This approach will not be
used elsewhere in this paper.
A regular integral polytope $\Delta \subset \R^n$ 
determines a smooth K\"ahler toric variety $(M,\omega)$,
and geometric quantization gives rise to a virtual representation $Q(M)$ 
of the torus $T^n$.
The dimension $\dim Q(M)$ of this quantization
is equal to the number of lattice points in $\Delta$. 
(This result (see \cite[Corollary 2.23]{Oda}) 
is an expression of the ``quantization commutes with reduction"
principle in symplectic geometry \cite{GS}.
According to this principle, $\dim Q(M)^c = \dim Q(M_c)$
for each lattice point $c \in \Z^n \subset {\rm Lie}(T^n)^*$,
where $Q(M)^c$ is the subspace of $Q(M)$ on which
$T^n$ acts through the character given by $c$,
and where $M_c$ is the reduced space of $M$ at $c$.  Because $M$ 
is a toric variety, $M_c$ is a point if $c \in \Delta$ 
and is empty otherwise.)
On the other hand, by the Hirzebruch-Atiyah-Singer generalization
of the classical Riemann-Roch formula, we have
$\dim Q(M) = \int_M \exp(c_1(L)) \Td(TM)$,
where $c_1(L) = [\omega]$ is the Chern class of the 
pre-quantization line bundle
and $\Td(TM)$ is the Todd class of the tangent bundle.
Expressing $M$ as a reduction of a linear torus action
on $\C^d$ (where $d$ is the number of facets of $\Delta$),
the tangent bundle stably splits into line bundles
$L_1,\ldots,L_d$,
and the above integral is obtained by applying
the Khovanskii-Pukhlikov differential operator $\prod \Td(\deldel{h_i})$
to the integral $\int_M \exp(\omega + \sum h_i c_1(L_i))$.
The Duistermaat-Heckman theorem on the variation
of reduced symplectic structures implies that this integral
is equal to the volume of the polytope $\Delta(h)$
that is obtained from $\Delta$ by shifting the $i$th facet 
by a distance $h_i$, for $i=1,\ldots,d$.
Hence, the number of lattice points in $\Delta$
is obtained by applying the Khovanskii-Pukhlikov operator
to the volume of $\Delta(h)$.

The Euler Maclaurin formulas due to Khovanskii-Pukhlikov,
Cappell-Shaneson, Guillemin, and Brion-Vergne are all
exact formulas, valid for sums of exponential or polynomial functions.
Cappell and Shaneson \cite{CS:private} have also investigated the
problem of deriving an Euler Maclaurin formula with remainder.
In a previous paper \cite{PNAS}, we stated and proved an Euler Maclaurin
formula with remainder for the sum of the values of an arbitrary smooth
function on the lattice points in a {\em regular}
polytope, and adumbrated a generalization to the case of
simple integral polytopes.  The purpose of this paper is to state
and prove an Euler Maclaurin formula with remainder
for simple lattice polytopes (Theorem \ref{main}).  The key
ingredients in the proof of this theorem, as in \cite{PNAS},
are a variant of the Euler Maclaurin formula in one dimension, 
given in Proposition \ref{twistedemray}, which
by iteration gives a formula for \emph{orthants}, 
along with a combinatorial result, given 
in Proposition \ref{weightedlawrence}, which shows how the sum 
of the values of
a function over the lattice points in a polytope can be decomposed into sums
over such orthants.

Our Euler Maclaurin formula with remainder is stated in Theorem \ref{main}
for functions of compact support. 
In Section \ref{sec:estimates} we show how to extend it to symbols 
(in the sense of Hormander, see e.g.\ \cite{horm}).
This is the content of Theorem \ref{syms}. 
As a corollary, we deduce an exact Euler Maclaurin formula for polynomials.

The early references to the Euler Maclaurin formula are Euler
\cite{Eu} and Maclaurin \cite{Ma}.  
\label{Poisson-page}
Apparently, Poisson
\cite{Poi} was the first to give a remainder formula.
See also \cite{H}.

% -------------------------------------------------------------------------
\section{Weighted sums in one dimension}
% -------------------------------------------------------------------------
\labell{sec:EM dim=1}

\subsection*{Exact Euler Maclaurin}

Here is a brief proof of the exact Euler Maclaurin formula
\eqref{EM exact weighted}; cf.~\cite{BV}.
First, we prove this formula
when $f(x)$ is an exponential function: $f(x) = e^{\lambda x}$
with $|\lambda| < 2 \pi.$  The formula then becomes
\begin{equation} \labell{EM for exponential}
\sum_{[a,b]}{'} e^{\lambda x}
   =   \left. \L(\deldel{h_1}) \L(\deldel{h_2}) \int_{a-h_1}^{b+h_2}
       e^{\lambda x} dx \right|_{h_1 = h_2 = 0} .
\end{equation}
An explicit computation,
which uses the facts that the constant term
in the formal power series $\L(S)$ is one and that $\L(-S) = \L(S)$, 
shows that 
\begin{equation} \labell{equation for LN}
 \left. \Ltk(\deldel{h_1}) \Ltk(\deldel{h_2}) \int_{a-h_1}^{b+h_2}
       e^{\lambda x} dx \right|_{h_1 = h_2 = 0}
 = \Ltk(\lambda) \int_a^b e^{\lambda x} dx 
\end{equation}
where $\Ltk(S)$ is the truncation of the power series $\L(S)$
at the even integer $2k$.
The radius of convergence of the power series $\L(\lambda)$ is $2\pi$
because the zeros of $\tanh(\lambda/2)$ that are nearest
to the origin are at $\pm 2 \pi i$.  Hence, for $|\lambda| < 2\pi$,
the expression in \eqref{equation for LN} converges as $k\to \infty$ to 
\begin{equation} \labell{LL1}
\L(\lambda) \int_a^b e^{\lambda x} dx
 = \frac{ \lambda/2 }{ \tanh(\lambda/2) } 
   \left( \frac{ e^{\lambda b} }{\lambda} 
        - \frac{ e^{\lambda a} }{\lambda} \right)
 = \half \left( e^{\lambda/2} + e^{-\lambda/2} \right)
   \frac{ e^{\lambda b} - e^{ \lambda a } }{ e^{\lambda/2} - e^{-\lambda/2} },
\end{equation} 
which is equal to the left hand side of \eqref{EM for exponential}
by the formula for a geometric sum.

Moreover, this convergence is uniform on any closed sub-disk
of $|\lambda| < 2\pi$, because $\L(\lambda)$ is a power series
and $\int_a^b e^{\lambda x} dx$ is bounded away from $0$ and from $\infty$.
Recall that differentiation commutes with uniform limits of holomorphic
functions (as a consequence of the Cauchy formula).
It follows that the derivative $\deldel{\lambda}$
commutes with the infinite order differential operators $\L(\deldel{h_i})$
on the right hand side of \eqref{EM for exponential}.
Comparing the Taylor coefficients of $\lambda^n$ on the left and right
hand sides of \eqref{EM for exponential}, we get a similar formula
for $f(x) = x^n$, and hence for all polynomials.

\begin{Remark} \labell{exp to pol}
In higher dimensions, the problem of obtaining a formula for polynomials
from a formula for exponentials is addressed in \cite{KP2}.
\end{Remark}

\subsection*{Weighted sums and ordinary sums}

The Todd function is defined by
$$ \Td(S) = \frac{S}{1 - e^{-S}}.$$
The Bernoulli numbers are the coefficients $b_k$ in
$$ \Td(-S) = \frac{S}{e^S - 1} 
           = 1 + \sum_{k \geq 1} \frac{1}{k!} b_k S^k . $$

(We are following the conventions in Bourbaki \cite{Bo}.)
Since
\begin{equation*}
\begin{aligned}
 \L(S) &= \left( S/2 \right) \frac{e^{S/2} + e^{-S/2}}{e^{S/2} - e^{-S/2}}
       = \half \left( \frac{S}{1 - e^{-S}} + \frac{S}{e^S - 1} \right)\\
       &= \half \left( \Td(S) + \Td(-S) \right) ,
\end{aligned}
\end{equation*}
the coefficients $b_{2k}$ in the power series expansion \eqref{LS}
are the even Bernoulli numbers.
Since
$$ \Td(S) - \Td(-S) = \frac{S}{1 - e^{-S}} - \frac{S}{e^S - 1} = S,$$
we have $b_{2k+1} = 0$ for all $k \geq 1$, so the only difference 
between $\Td(S)$ and $\L(S)$
is the absence of the linear term in $\L(S)$. 
Replacing $\L(\cdot)$ by $\Td(\cdot)$ on the right hand side
of the Euler Maclaurin formula \eqref{EM exact weighted}
results in a formula for the ordinary sum of the values of $f$
over the integers in $[a,b]$.
However, we work with the formula for the weighted sum \eqref{wted sum},
because this formula avoids ``boundary effects" which occur 
in the formulas for the ordinary sum.

\subsection*{Euler Maclaurin with remainder}

As before,
we will denote the truncation of the power series $\L(S)$
at the even integer $2k$ by $\Ltk(S)$. Then
\begin{equation} \labell{Ltk}
 \Ltk(\deldel{h}) 
   = 1 + \frac{1}{2!} b_2 \frac{\del^2}{\del h^2}
       + \ldots 
       + \frac{1}{(2k)!} b_{2k} \frac{\del^{2k}}{\del h^{2k}}
\end{equation}
is a differential operator with constant coefficients involving only even
order derivatives.  In particular, 
if $g(h)$ is a function with $2k$ continuous derivatives, then
\begin{equation} \labell{plusminus}
\left.\Ltk(\deldel{h})(g(h))\right|_{h=0} 
  = \left.\Ltk(\deldel{h})(g(-h))\right|_{h=0}.
\end{equation}

The classical Euler Maclaurin summation formula with remainder  
can be formulated in the following way.
\begin{Proposition}[Euler Maclaurin with remainder for intervals]
\labell{EMintervalgrm}
Let $f(x)$ be a function with $m \geq 1$ continuous derivatives
and let $k = \lfloor m/2 \rfloor$.
Then
\begin{multline} \labell{emr1d}
  {\sum_{[a,b]} }' f = \left. \Ltk(\deldel{h_1}) \Ltk(\deldel{h_2})
  \int_{a-h_1}^{b+h_2} f(x) dx
  \right|_{h_1 = h_2 = 0} 
 + (-1)^{m-1} \int_a^b P_m(x) f^{(m)}(x) dx 
\end{multline}
with
\begin{equation} \labell{P and B}
 P_m(x) = \frac{B_m(\{ x \})}{m!},
\end{equation}
where $B_m(x)$ is the $m$th Bernoulli polynomial (see below)
and where $\{ x \} = x - \lfloor x \rfloor$ 
is the fractional part of $x$.
Moreover, the function $P_m(x)$ is given by 
\begin{equation} \labell{P2k} 
 P_{2k}(x) = (-1)^{k-1} \sum_{n=1}^\infty \frac{2\cos2\pi nx}{(2\pi n)^{2k}} 
\end{equation}
if $m=2k$ is even and by
\begin{equation} \labell{P2kp1}
 P_{2k+1}(x) =
   (-1)^{k-1} \sum_{n=1}^\infty \frac{2\sin2\pi n x}{(2\pi n)^{2k+1}}
\end{equation}
if $m=2k+1$ is odd.
\end{Proposition}

Up to minor changes in notation, this result is formula \textbf{298} 
in \cite{Kn}.  

Note that if $f$ is a polynomial 
then \eqref{emr1d} becomes an exact formula when $m$ is greater
than the degree of $f$. 

\medskip

Let us recall the proof of Proposition \ref{EMintervalgrm}.
Consider the difference between integration and summation 
as the difference between two distributions:
$$ P_0 (x) := 1 - \sum_{k \in \Z} \delta(x-k). $$
Integrating this, with the choice of constant of integration 
such that the integral from $0$ to $1$ 
of the result vanishes, gives a distribution $P_1$ which
is given by the function
\begin{equation} \labell{P1x}
 P_1(x) = x - \lfloor x \rfloor - \half
\end{equation}

\noindent at non-integral points.
This is the famous zigzag function studied by Gibbs \cite{Gibbs}.  
\label{Gibbs-page}
Notice that $P_1$ is 1-periodic and odd.

If $f$ is a continuously differentiable function on $[0,1]$, we get
$$
\int_0^1 P_1(x) f'(x) dx 
  = \left. \phantom{\int} P_1(x) f(x) \right|_0^1 - \int_0^1 f(x) dx 
  = \half f(0) + \half f(1) - \int_0^1 f(x) dx,
$$
so 
$$\half f(0) + \half f(1) = \int_0^1 f(x) dx + \int_0^1 P_1(x) f'(x) dx.$$
Summing the corresponding expressions over the intervals
$[a,a+1]$, $[a+1,a+2]$, $\ldots$, $[b-1,b]$, 
where $a <b$ are integers, gives
\begin{equation} \labell{sum over intervals}
 \sum_{[a,b]}{'} f = \int_a^b f(x) dx + \int_a^b P_1(x) f'(x) dx. 
\end{equation}
This is the first of a series of expressions relating the sum to an integral. 

We now take successive anti-derivatives:
for $m \geq 2$,
define $P_m(x)$ inductively by the conditions that
$\frac{d}{dx} P_m(x) = P_{m-1}(x)$ and $\int_0^1 P_m(x) dx = 0$.
Then $P_m(x)$ is 1-periodic and satisfies $P_m(-x) = (-1)^m P_m(x)$.
Also, $P_m(x)$ is continuously differentiable $m-2$ times;
in particular, it is continuous if $m \geq 2$.
Starting from \eqref{sum over intervals}, we integrate by parts:
\begin{eqnarray*}
\sum_{[a,b]}{'} f
 & = & \int_a^b f(x) dx + \left. P_2 f' \right|_a^b - \int_a^b P_2(x) f''(x) dx \\
 & = & \int_a^b f(x) dx + \left. P_2 f' \right|_a^b 
                        - \left. P_3 f'' \right|_a^b + \int_a^b P_3(x) f^{(3)}(x) dx\\
 & = & \int_a^b f(x) dx + \left. P_2 f' \right|_a^b - \left. P_3 f'' \right|_a^b
                        + \left. P_4 f^{(3)} \right|_a^b - \int_a^b P_4(x) f^{(4)}(x) dx \\
 & \vdots & 
\end{eqnarray*}
Noting that $P_n(a) = P_n(b) = P_n(0)$ and $P_{2n+1}(0) = 0$
(because $P_{2n+1}$ is odd), we get, setting $k = \lfloor m/2 \rfloor$,
\begin{multline} \labell{dog}
 \sum_{[a,b]}{'} f  =  
       \int_a^b f(x) dx + P_2(0) f' |_a^b + P_4(0) f^{(3)}|_a^b + \ldots
       + P_{2k}(0) f^{(2k-1)}|_a^b \\
 + (-1)^{m-1} \int_a^b P_m(x) f^{(m)}(x) dx.$$
\end{multline}

Consider the polynomial
\begin{equation} \labell{L2k[]}
 L^{[2k]}(S) := 1 + P_2(0) S^2 + P_4(0) S^4 + \ldots + P_{2k}(0) S^{2k}.
\end{equation}
From \eqref{dog} we get the remainder formula 
\begin{multline} \labell{emr1d[]}
   \sum_{[a,b]}{'} f = 
   L^{[2k]} (\deldel{h_1}) L^{[2k]} (\deldel{h_2}) 
   \left. \int_{a \pm h_1}^{b \pm h_2} f(x) dx \right|_{h_1 = h_2 = 0}
   + (-1)^{m-1} \int_a^b P_{m}(x) f^{(m)}(x) dx 
\end{multline}
for a function $f$ of type $C^m$, where $k = \lfloor m/2 \rfloor$.

This formula becomes exact when $f$ is a polynomial
and when $m$ is sufficiently large.
We therefore have, by comparison with Equation \eqref{EM exact weighted},
\begin{equation} \labell{L[]=L}
L^{[2k]}= \Ltk.
\end{equation}
This and \eqref{emr1d[]} give \eqref{emr1d}.

It remains to derive the expressions \eqref{P and B},
\eqref{P2k}, and \eqref{P2kp1} for the functions $P_{m}(x)$.

The Bernoulli polynomials $B_m(x)$ are characterized by the properties
$$ B_0(x) =1 \ , \ \ \frac{d}{dx} B_m(x) = m B_{m-1}(x) \ , \ \ 
\text{and} \ \int_0^1 B_m(x) dx = 0 \ \ \text{ for } \ m \geq 1.$$
In particular, 
$$ B_1(x) = x - \half = P_1(x) \quad \text{ for all } \ 0 < x < 1. $$
(See \eqref{P1x}.)
Integrating, we get that $P_m(x) = B_m(x) / m!$ for all $0 < x < 1$,
and since $P_m$ is 1-periodic, we get \eqref{P and B}.

The Poisson summation formula says that, as distributions,
$$ \sum_{k \in \Z} \delta(x-k) = \sum_{n \in \Z} e^{2\pi i n x},$$
so
$$ P_0(x) = 1 - \sum_{k \in \Z} \delta(x-k) 
= -2 \sum_{n=1}^\infty \cos 2\pi n x,$$
as distributions.
Integrating this, with the choice of constant of integration 
such that the integral from $0$ to $1$ of the result vanishes, 
gives the expressions \eqref{P2k} and \eqref{P2kp1}.

The idea of using these Fourier expansions goes back to Wirtinger \cite{Wi}.
See Knopp \cite{Kn}, pp.~521--524.
\label{Wi-page}

\begin{Remark}
From \eqref{L2k[]} and \eqref{P2k} we see that 
the coefficients of the polynomials $L^{[2k]}(S)$ are 
$$\displaystyle { P_{2k}(0) = 
   (-1)^{k-1} \frac{2}{(2\pi)^{2k}} \sum_{n \geq 1} \frac{1}{n^{2k}} }
 = (-1)^{k-1} \frac{2}{(2\pi)^{2k}} \zeta(2k),$$
where $\zeta(\cdot)$ is the Riemann zeta function.
Comparing this with the coefficients of the polynomials $L^{2k}(S)$,
we get, from \eqref{L[]=L} and \eqref{Ltk}, that
$$\zeta(2k) = (-1)^{k-1} \frac{(2\pi)^{2k}}{2} \frac{1}{(2k)!} b_{2k},$$
reproducing Euler's famous evaluation 
of the Bernoulli numbers in terms of the Riemann zeta function.
\label{zeta-page}
\end{Remark}

Proposition \ref{EMintervalgrm}, when applied to a $C^{m}$ function $f(x)$ 
of compact support, implies a similar formula for an infinite ray:
\begin{multline} \labell{emr1dray}
 \half f(a) + f(a+1) + f(a+2) + \ldots \\
 = \left. \Ltk(\deldel{h}) \int_{a\pm h}^{\infty} f(x) dx
 \right|_{h=0} 
 + (-1)^{m-1} \int_a^\infty P_{m}(x) f^{(m)}(x) dx$$
\end{multline}
for either choice of $\pm$, where $k = \lfloor m/2 \rfloor$.
Indeed, we need only choose $b$ so large that the support of $f(x)$ 
is contained in the set $\{ x < b \}$,
and then apply Proposition \ref{EMintervalgrm}.
Conversely, if we know (\ref{emr1dray}) for
functions of compact support, then we can conclude (\ref{emr1d}). 
Indeed, multiply $f$ by a smooth function of compact support which is
equal to one in a neighborhood of $[a,b]$ and observe that
\begin{equation}
\labell{1dweighedLawrence}
{\sum_{[a,b]} }'={\sum_{[a,\infty)} }'-{\sum_{[b,\infty)} }'
\end{equation}
when applied to a function of compact support.

\subsection*{Twisted Euler Maclaurin with remainder for a ray.}

In extending the Euler Maclaurin formula to higher dimensions, 
we will need an expression for the ``twisted weighted sum"
$$ \half f(0) + \sum_{n=1}^\infty \lambda^n f(n) $$ 
when $\lambda \neq 1$ is a root of unity, say, of order $N$,
in terms of the integrals of $f$.

Let
\begin{equation} \labell{def:Q0} 
 Q_{0,\lambda}(x) = - \sum_{n \in \Z} \lambda^n \delta(x-n). 
\end{equation}
This is an $N$-periodic distribution since $\lambda^N = 1$.
We will take successive anti-derivatives of this distribution,
with the constants chosen so that the integrals from $0$ to $N$ vanish:
\begin{equation} \labell{def:Qm}
 \frac{d}{dx} Q_{m,\lambda} (x) = Q_{m-1,\lambda} (x)
\quad \text{and} \quad \int_0^N Q_{m,\lambda}(x) dx = 0 . 
\end{equation}
With this choice of constants, $Q_{m,\lambda}(x)$ is $N$-periodic 
for each $m$.

Let us now look more  closely at $Q_{1,\lambda}(x)$. We have
$$ \frac{d}{dx} {\bf 1}_{[n,n+1]}(x) = \delta(x-n)-\delta(x-(n+1)),$$
so
$$ \frac{d}{dx} \sum_{n \in \Z} \lambda^n{\bf 1}_{[n,n+1]} (x) =
   \sum_{n \in \Z} (\lambda^n - \lambda^{n-1}) \delta(x-n)
   = \frac{1-\lambda}{\lambda} Q_{0,\lambda} (x).$$
Thus,
$$ Q_{1,\lambda}(x) = \frac{\lambda}{1-\lambda} 
   \sum_{n \in \Z} \lambda^n {\bf 1}_{[n,n+1]} (x),
$$
because  the right hand side is periodic of period $N$ 
and its integral over $[0,N]$ vanishes since 
$1 + \lambda + \lambda^2 + \cdots + \lambda^{N-1} = 0$. 
If $f$ is a continuously differentiable function of compact support,
we have
\begin{multline*}
 \int_0^\infty Q_{1,\lambda}(x) f'(x) dx = 
   \frac{\lambda}{1-\lambda} \left. \sum_{n=0}^\infty \lambda^n f 
   \right|_n^{n+1} 
  = - \frac{\lambda}{1-\lambda} f(0) + \lambda f(1) + \lambda^2 f(2) 
    + \cdots 
\end{multline*}
so
\begin{equation} \labell{firsttwistedemray}
\half f(0) + \sum_{n\geq 1} \lambda^n f(n) = 
 \left( \half + \frac{\lambda}{1-\lambda} \right) f(0)
 + \int_0^\infty Q_{1,\lambda}(x) f'(x) dx.
\end{equation}

Successively applying integration by parts to \eqref{firsttwistedemray}, 
we get
\begin{multline*}
 \half f(0) + \lambda f(1) + \lambda^2 f(2) + \ldots \\
    = \left( \half + \frac{\lambda}{1 - \lambda} \right) f(0)
      - Q_{2,\lambda}(0) f'(0) - \int_0^\infty Q_{2,\lambda}(x) f''(x) dx \\
    = \left( \half + \frac{\lambda}{1 - \lambda} \right) f(0)
      - Q_{2,\lambda}(0) f'(0) + Q_{3,\lambda}(0) f''(0)
      + \int_0^\infty Q_{3,\lambda}(x) f^{(3)}(x) dx \\
     \vdots \\
    \begin{aligned}
   = \left( \half + \frac{\lambda}{1 - \lambda} \right) f(0)
      - Q_{2,\lambda}(0) f'(0) + \ldots 
      + (-1)^{k-1} Q_{k,\lambda}(0) f^{(k-1)}(0) \\
      \quad + (-1)^{k-1} \int_0^\infty Q_{k,\lambda}(x) f^{(k)}(x) dx.
\end{aligned}
\end{multline*}
Since 
$$ (-1)^{m-1} f^{(m-1)}(0) = \left. 
  \left( \deldel{h} \right)^{m} \int_{-h}^\infty f(x) dx \right|_{h=0},$$
we have thus proved this ``twisted Euler Maclaurin formula for a ray":

\begin{Proposition}\labell{twistedemray}
Let
$${\bfM}^{k,\lambda}(S) =
\left( \frac12 + \frac{\lambda}{1-\lambda} \right) S
 + Q_{2,\lambda}(0) S^2
 + Q_{3,\lambda}(0) S^3
 + \cdots 
 + Q_{k,\lambda}(0) S^k,$$
for a root of unity $\lambda \neq 1$.  Then
\begin{multline} \labell{twistedemrayeq}
 \half f(0) + \lambda f(1) + \lambda^2 f(2) + \cdots = 
  \left. 
  {\bfM}^{k,\lambda}(\deldel{h}) \int_{-h}^\infty f(x) dx \right|_{h=0} \\
  + (-1)^{k-1}\int_0^\infty Q_{k,\lambda}(x) f^{(k)}(x) dx
\end{multline}
if $f \in C_c^{k}(\R)$.
\end{Proposition}

\begin{Remark}
Another twisted Euler Maclaurin formula for a ray appeared in \cite{T}. 
\end{Remark}

We will now show that the polynomials that appear in Proposition 
\ref{twistedemray} have the following symmetry property:
\begin{equation}
\labell{Dandlambdainverse}
\bfM^{m,\lambda^{-1}}(S)= \bfM^{m,\lambda}(-S).
\end{equation}

The linear term transforms according to \eqref{Dandlambdainverse} 
because
$$ \half + \frac{\lambda\inv}{1 - \lambda\inv}
   = \half - \left( 1 + \frac{\lambda}{1 - \lambda} \right)
 = - \left( \half + \frac{\lambda}{1 - \lambda} \right).$$
%
%$$\half + \frac{\lambda}{1-\lambda} 
 %+ \half + \frac{\lambda^{-1}}{1-\lambda^{-1}}
 %= 1 + \frac{\lambda}{1-\lambda} + \frac1{\lambda-1} = 0.$$
For the other terms to transform correctly, we need to check that
\begin{equation} \labell{lambdaolambdainv}
Q_{m,\lambda^{-1}}(0) = (-1)^m Q_{m,\lambda}(0) 
  \quad \text{for all} \quad m \geq 2.
\end{equation}

As in the non-twisted formula, 
it is convenient to work with the Fourier expansions of the
$Q_{m,\lambda}$'s.  Let
$$ \lambda = e^{2\pi i \frac{j}{N}}.$$
We are assuming that $\lambda \neq 1$, and so $j \not \equiv 0 \mod N$.
Then we can write the distribution $Q_{0,\lambda}(x)$ as 
$$ Q_{0,\lambda}(x) = - \sum_{n \in \Z} \lambda^n \delta(x-n)
   = - e^{ (2 \pi i \frac{j}{N}) x } \sum_{n \in \Z} \delta(x-n). $$
Writing
$$ \sum_{n \in \Z} \delta(x-n) = \sum_{r \in \Z} e^{2\pi i rx },$$
we see that the Fourier series of $Q_{0,\lambda}(x)$ is
$$ Q_{0,\lambda}(x) = - \sum_{r \in \Z} e^{ 2 \pi i (\frac{j}{N} +r) x}. $$
The indefinite integral, chosen so that the integral over $[0,N]$ 
vanishes, is obtained by dividing each Fourier summand by the coefficient 
of the exponent, so we have the following Fourier series:
\begin{equation}\labell{fourierq}
 Q_{m,\lambda}(x) = 
   - \sum_{r \in \Z} \frac {e^{2\pi i(\frac{j}{N}+r)x}} 
                {\left(2\pi i(\frac{j}{N}+r)\right)^m}.
\end{equation}
Setting $x=0$, replacing $j$ by $-j$, and replacing $r$ by $-r$ 
in the sum, pulls out a factor of $(-1)^m$, 
which gives the required equation \eqref{lambdaolambdainv},
and which implies \eqref{Dandlambdainverse}.

\begin{Remark} \labell{M vs L}
For $\lambda = 1$, if we define
$$ \bfM^{k,1}(S) = \L^{2 \lfloor k/2 \rfloor}(S) \quad \text{and} \quad
   Q_{k,1} = P_k,$$
then \eqref{twistedemrayeq} boils down to \eqref{emr1dray}.
So, with this notation, \eqref{twistedemrayeq} also holds for 
$\lambda = 1$.  
Notice that if $\lambda \neq 1$ then $M^{k,\lambda}(S)$ is a multiple
of  $S$, and that if $\lambda = 1$ then $M^{k,\lambda}(S) = 1 + $
a multiple of $S$.
Finally, the symmetry property \eqref{Dandlambdainverse} continues to hold
for $\lambda = 1$ because the polynomials $\L^{2k}$ are symmetric.
\end{Remark}

% ----------------------------------------------------------------
\section{The polar decomposition of a simple polytope}
% ----------------------------------------------------------------
\labell{sec:decompose}

In this section we decompose a polytope into an ``alternating sum
of polyhedral cones". See \cite{V} or \cite{lawrence}.
We will give a ``weighted version" of this decomposition.

Our polytopes are always compact and convex.
A compact convex polytope $\Delta$ in $\R^n$ is a compact set
which can be obtained as the intersection of finitely many
half-spaces, say,
\begin{equation}\labell{halfspaces}
 \Delta = H_1 \cap \ldots \cap H_d.
\end{equation}
We assume that \eqref{halfspaces} is an intersection
with the smallest possible $d$, so that the $H_i$'s are uniquely 
determined up to permutation.  We order them arbitrarily.  
The \emph{facets} (codimension one faces) of $\Delta$ are
$$ \sigma_i = \Delta \cap \del H_i \quad , \quad i = 1, \ldots, d .$$

Alternatively, a compact convex polytope is the convex hull
of a finite set of points in $\R^n$.  If we take this set to be minimal,
it is uniquely determined, and its elements are the \emph{vertices}
of $\Delta$.

For each vertex $v$ of $\Delta$, let $I_v \subset \{ 1, \ldots, d \}$
encode the set of facets that contain $v$, so that 
$$ i \in I_v \quad \text{if and only if} \quad v \in \sigma_i .$$

Assume that $\Delta$ is \emph{simple}, so that each vertex is
the intersection of exactly $n$ facets.
For each $i \in I_v$, there exists a unique edge at $v$ which does not 
belong to the facet $\sigma_i$; choose any vector $\alpha_{i,v}$ 
in the direction of this edge.
(At the moment, these ``edge vectors" are only determined
up to a positive scalar. Later, when discussing integral polytopes, 
we will make a specific choice of the edge vectors.)

The \emph{tangent cone} to $\Delta$ at $v$ is
$$
 \bfC_v =
    \{ v + r (x-v) \ \mid \ r \geq 0 \, , \, x \in \Delta \} 
 = v + \sum_{j \in I_v} \R_{\geq 0} \alpha_{j,v}.
$$
We will ``polarize" these tangent cones by flipping some of their edges
so that they all ``point in the same direction".
This direction is specified by the choice of a ``polarizing vector":
a vector
$\xi \in {\R^n}^*$, such that $\langle \xi , \alpha_{j,v} \rangle$
is non-zero for all $v$ and $j$.  With this choice, we define
the \emph{polarized edge vectors} to be
\begin{equation} \labell{flipalpha}
  \valpha^{\sharp}_{i,v} = \begin{cases}
 \valpha_{i,v} & \text{ if } \l< \xi , \valpha_{i,v} \r> < 0, \\
 - \valpha_{i,v} & \text{ if } \l< \xi , \valpha_{i,v} \r> > 0 ,
 \end{cases}
\end{equation}
and the \emph{polarized tangent cone} to be
\begin{equation} \labell{pol tan cone}
 \bfC_v^\sharp = v + \sum_{j \in I_v} \R_{\geq 0} \alpha_{j,v}^\sharp. 
\end{equation}

We define the ``weighted characteristic function", 
\begin{equation} \labell{def:wted char}
{\bf 1}^w_\Delta(x),
\end{equation}
to be the function on $\R^n$ that takes the value $0$ on the exterior 
of $\Delta$, the value $1$ on the interior of $\Delta$, and 
the value $1/2^k$ on the relative interior of a codimension $k$ face 
of $\Delta$.
So, for example, for an interval $[a,b]$ on the line, the function
${\bf 1}^w_{[a,b]}(x)$ assigns the value
one to points $a<x<b$, zero to points outside the interval, 
and $\frac12$ to the points $a$ and $b$.
We use a similar definition for a polyhedral cone.

\begin{Proposition}[\protect{Weighted polar decomposition 
of a simple polytope}]
\labell{weightedlawrence}
Let $\Delta$ be a simple polytope.
For any choice of polarizing vector, we have
\begin{equation} \labell{wted L}
 {\bf 1}^w_\Delta (x) = \sum_v(-1)^{\#v} {\bf 1}^w_{\bfC_v^\sharp}(x),
\end{equation}
where the sum is over the vertices $v$ of $\Delta$,
where $\bfC_v^\sharp$ is the polarized tangent cone,
where ${\bf 1}^w_\Delta (x)$
and ${\bf 1}^w_{\bfC_v^\sharp}(x)$ are the weighted characteristic
functions, and where $\#v$ denotes the number of 
edge vectors at $v$ whose signs were changed by
the polarizing process \eqref{flipalpha}.
\end{Proposition}

\begin{figure}
\setlength{\unitlength}{0.00004in}
\begin{center}
\psfig{figure=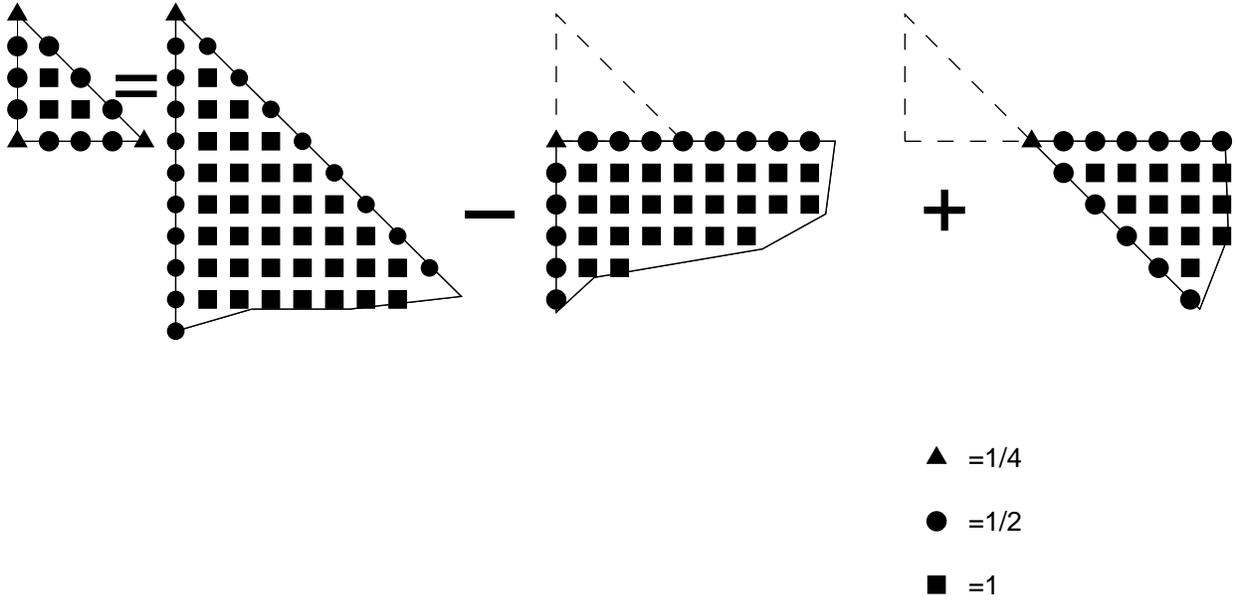,width=\textwidth}
\end{center}
\caption{Polar decomposition of a triangle}
\label{fig:Lawrence}
\end{figure}

This theorem is illustrated for the case of a triangle in Figure 
\ref{fig:Lawrence}.

\vskip.2in
\noindent
\begin{proof}

We will prove this equality in two steps.
First, for each $x$ we will find a polarization
such that the equality \eqref{wted L} holds. Second, 
we will show that the right hand side of the equality 
\eqref{wted L} is independent of the choice of polarization.

Suppose that $x \not \in \Delta$.
Let $\xi$ be a polarizing vector such that 
$\left< \xi , x \right> > \left< \xi , y \right> $
for all $y \in \Delta$. Then all the polarized cones
``point away from $x$" and do not contain $x$.
Formula (\ref{wted L}) for the polarizing vector $\xi$,
when evaluated at $x$, states that $0=0$.

For every polarization there exists exactly one vertex $v$
for which $\bfC_v^\sharp = \bfC_v$, namely, the vertex $v$
such that $\langle \xi , v \rangle$ is maximal. Conversely, for every
vertex $v$ there exists a polarization such that $\bfC_v^\sharp = \bfC_v$.

Suppose that $x \in \Delta$ is in the relative interior
of a face $F$ of codimension $k$.  Let $v$ be any vertex of $F$. 
Let $\xi$ be a polarizing vector such that $\bfC_v^\sharp = \bfC_v$. 
Let $F_v$ be the codimension $k$ face of $\bfC_v$ that contains $F$.
Then $ {\mathbf 1}_\Delta^w(x) 
  = {\mathbf 1}_{C_v^\sharp}^w (x) = \left( \frac{1}{2} \right)^k$.
For each other vertex, $v'$, the cone $C_{v'}^\sharp$ is disjoint
from the relative interior of $F$, so ${\textbf 1}_{C_v^\sharp}^w(x) = 0$. 
Formula (\ref{wted L}) for this polarization, when evaluated at $x$, 
states that $\left(\half\right)^k = \left(\half\right)^k$.

To keep track of the different possible choices of polarization, we let
$$ E_1,\ldots,E_N $$
denote all the different codimension one subspaces of ${\R^n}^*$
that are equal to
$$ \alpha_{j,v}^\perp = 
   \{ \eta \in {\R^n}^* \mid \langle \eta , \alpha_{j,v} \rangle =0 \} $$
for some $j$ and $v$.
(For instance, if no two edges of $\Delta$ are parallel,
then the number $N$ of such hyperplanes
is equal to the number of edges of $\Delta$.)
A vector $\xi$ can be taken to be a ``polarizing vector"
if and only if it does not belong to any $E_j$.
The ``polarized cones" $\bfC_v^\sharp$ only depend on the
connected component of the complement 
\begin{equation} \labell{tDelta}
 {\R^n}^* \ssminus (E_1 \cup \ldots \cup E_N) 
\end{equation}
in which $\xi$ lies. 
Any two polarizing vectors can be connected by a path $\xi_t$
in ${\R^n}^*$ which crosses the ``walls" $E_j$ one at a time.
We finish by showing that the right hand side of formula (\ref{wted L})
does not change when the polarizing vector $\xi_t$ crosses
a single wall, $E_k$.

As $\xi_t$ crosses the wall $E_k$,
the sign of the pairing $\langle \xi_t , \alpha_{j,v} \rangle$ flips
exactly if $E_k = \alpha_{j,v}^\perp$. 
For each vertex $v$, denote by $\bfS_v(x)$ and $\bfS_v'(x)$
its contributions to the right hand side of formula
\eqref{wted L} before and after $\xi_t$ crossed the wall.
The vertices for which these contributions differ are exactly those
vertices that lie on edges $e$ of $\Delta$ which are perpendicular
to $E_k$.  They come in pairs because each edge has two endpoints.

Let us concentrate on one such an edge, $e$, with endpoints, say,
$u$ and $v$. Let $\alpha_e$ denote an edge vector at $v$ that points 
from $v$ to $u$ along $e$.
Suppose that the pairing $\langle \xi_t , \alpha_e \rangle $ 
flips its sign from negative to positive as $\xi_t$ crosses the wall;
(otherwise we switch the roles of $v$ and $u$).
The ``polarized tangent cones" to $\Delta$ at $v$, before and after
$\xi_t$ crosses the wall, are
\begin{equation*} \begin{aligned}
 \bfC_v^\sharp &= v + \sum_{j \in I_e} \R_{\geq 0} \alpha_{j,v}^\sharp
                + \R_{\geq 0} \alpha_e 
\quad \text{and} \\
 (\bfC_v^\sharp)' & = v + \sum_{j \in I_e} \R_{\geq 0} \alpha_{j,v}^\sharp
                - \R_{\geq 0} \alpha_e
\end{aligned} \end{equation*}
where $I_e \subset \{ 1, \ldots, d\}$ encodes the facets that contain $e$.
(The $\alpha_{j,v}^\sharp$ are the same for the different $\xi_t$'s
because the pairings 
$ \langle \xi_t , \alpha_{j,v}^\sharp \rangle $ do not flip sign
when $\xi_t$ crosses the wall for $j \in I_e$.) 
The cones $\bfC_v^\sharp$ and $(\bfC_v^\sharp)'$ have a common facet
and their union is
$$ \bfC_e^\sharp := v + \sum_{j \in I_e} \R_{\geq 0} \alpha_{j,v}^\sharp 
                      + \R \alpha_e.$$
This union only depends on the edge $e$ and not on the endpoint $v$.
(This uses the assumption that the polytope $\Delta$ is simple
and follows from the fact that 
$\alpha_{j,u} \in \R_{\geq 0} \alpha_{j,v} + \R \alpha_e$.)

The contributions of $v$ to the right hand side of (\ref{wted L}) 
before and after
$\xi_t$ crosses the wall are
\begin{equation} \labell{bfSv}
 \bfS_v(x) = \varepsilon {\bf 1}^w_{\bfC_v^\sharp} \quad \text{and} \quad
\bfS_v'(x) = -\varepsilon {\bf 1}^w_{(\bfC_v^\sharp)'}
\end{equation}
where $\varepsilon \in \{ -1 , 1 \}$.
Their difference is plus/minus the weighted characteristic function 
of $\bfC_e^\sharp$:
\begin{equation} \labell{difference}
 \bfS_v(x) - \bfS_v'(x) = \varepsilon {\bf 1}^w_{\bfC_e^\sharp}
\end{equation}
The contributions of the other endpoint, $u$, have opposite signs 
than the respective contributions \eqref{bfSv} of $v$,
and their difference is minus/plus the characteristic function 
of $\bfC_e^\sharp$.  Hence, the differences $\bfS_v(x)-\bfS_v'(x)$
and $\bfS_u(x)-\bfS_u'(x)$, for the two endpoints $u$ and $v$ of $e$, 
sum to zero.
\end{proof} 

% ---------------------------------------------------------------------
\section{Euler Maclaurin formula with remainder for regular polytopes.}
% ---------------------------------------------------------------------
\labell{sec:EM regular}
We extend our notation \eqref{wted sum}
to a weighted sum over a simple integral polytope~$\Delta$:
$$ \sum_{\Delta \cap \Z^n}{'} f := \sum_{x \in \Delta \cap \Z^n} 
   {\bf 1}_\Delta^w(x) f(x), $$
where ${\bf 1}_\Delta^w(x)$ is the weighted characteristic function,
introduced in \eqref{def:wted char}, which is equal to $1/2^k$ 
when $x$ lies in the relative interior of a face of $\Delta$ 
of codimension $k$.
We use similar notation for the weighted sum of a compactly supported
function $f$ over a simple polyhedral cone.

With this notation, Proposition \ref{weightedlawrence}
gives the following decomposition. Let $\xi \in {\R^n}^*$
be a ``polarizing vector". Then 
\begin{equation} \labell{wheee}
 \sum_{\Delta \cap \Z^n}{'} f = \sum_v (-1)^{\# v}
 \sum_{\bfC_v^\sharp \cap \Z^n}{'} f 
\end{equation}
where we sum over the vertices $v$ of $\Delta$,
where $\bfC_v^\sharp$ is the polarized tangent cone to $\Delta$ at $v$
(see \eqref{pol tan cone}), and where $\# v$ is the number of edge vectors
at $v$ that are flipped by the polarization process \eqref{flipalpha}.

\bigskip

For the standard closed orthant 
$\bfO = \prod_{i=1}^n \R_{\nonneg}$ in $\R^n$, we have  
$${\sum_{\bfO \cap \Z^n}}' g
 = {\sum_{m_1 \in \Z_{\nonneg}}}' \cdots
   {\sum_{m_n \in \Z_{\nonneg}}}' g(m_1,\dots,m_n) , $$
for any function $g$ of compact support.

By performing $n$ iterations of \eqref{emr1dray}, we obtain
an Euler Maclaurin formula with remainder for the standard orthant:
Let $m$ be an integer and let $k = \lfloor m/2 \rfloor$.
If $g$ is a $C^{mn}$ function of compact support and $m\geq 1$, then
for any choice of the $\pm$'s we have
\begin{equation}   \labell{EMstandardorthant}
 {\sum_{\bfO\cap \Z^n}}'  g = \left. \prod_{i=1}^n
 \Ltk(\deldel{h_i}) \int \limits_{\bfO(\pm h_1,\dots,\pm h_n)} g(x) dx
 \right|_{h_1=\cdots=h_n=0} + R^{\st}_{m}(g),
\end{equation}
where
$$\bfO(h_1,\dots,h_n)=\{t \mid t_i + h_i \geq 0 \ \forall i \}$$ 
denotes the shifted orthant, and where the remainder term
is given by  
\begin{multline} \labell{standardremainder1} 
 R^{\st}_{m}(g) = \\
 \sum_{I\subsetneq\{1,\dots,n\}} (-1)^{(m-1)(n-|I|)} \prod_{i\in I}
 \Ltk(\deldel{h_i})
\left. 
 \int\limits_{\bfO(\pm h_1,\dots,\pm h_n)}
 \prod_{i\not\in I}P_{m}(x_i)\prod_{i\not\in I}
 \left(\frac{\partial}{\partial x_i}\right)^{m} g(x) dx_1 \cdots dx_n
 \right|_{h=0} .
\end{multline}
This remainder can also be expressed as a sum of integrals
over the orthant of bounded periodic functions times various 
partial derivatives of $f$ of order no less than $m$ and no more than $mn$. 
This fact follows from the formula
\begin{equation} \labell{formula with varphi}
\deldel{h_i} \int\limits_{\bfO(h_1, \ldots, h_n)}
  \varphi dx_1 \cdots dx_n = 
  - \int\limits_{\bfO(h_1, \ldots, h_n)}
    \frac{\del \varphi}{\del x_i} dx_1 \cdots dx_n.
\end{equation}
We can apply iterations of this formula
with $i \in I$ to the $I$'th summand of \eqref{standardremainder1}
because the non-smooth functions $P_m$ are only applied to the
variables $x_j$ for $j \not \in I$. We get
\begin{equation} \labell{standardremainder}
R^{\st}_{m}(g) = \sum_{I \subsetneq \{ 1, \ldots, n \} }
   (-1)^{(m-1)(n-|I|)} \int\limits_\bfO 
\prod_{i \in I} \Ltk (-\deldel{x_i})
\prod_{i \not \in I} P_m(x_i)
\prod_{i \not \in I} (\deldel{x_i})^m 
g(x) dx_1 \ldots dx_n.
\end{equation}

\bigskip

A regular integral orthant $\bfC$ is the image of the standard orthant $\bfO$
via an affine transformation
of the form
\begin{equation} \labell{affine}
(t_1,\ldots,t_n) \mapsto x = v + t_1 \alpha_1 + \ldots + t_n \alpha_n
\end{equation}
where $\alpha_1,\ldots,\alpha_n$ generate $\Z^n$
and where $v \in \Z^n$.
If $u_1,\ldots,u_n \in {\R^n}^*$ 
is the dual basis to $\alpha_1,\ldots,\alpha_n$
then the image of $\bfO(h_1,\ldots,h_n)$ under this transformation
is given by the inequalities
\begin{equation}\labell{uiss}
\langle u_i , x \rangle-\langle u_i , v \rangle+h_i\geq 0.
\end{equation}
We denote this expanded orthant by ${\bf C}(h)$. If $f$ is a $C^{mn}$
function of compact support and
$$ g(t_1,\ldots,t_n) = f(t_1 \alpha_1 + \ldots + t_n \alpha_n) $$
is its pullback under the transformation \eqref{affine}, then
$${\sum_{{\bfC}\cap \Z^n}}'f={\sum_{\bfO \cap \Z^n}}'g \quad \text{and} \quad
  \int\limits_{{\bfC}(h)}f (x) dx
 =\int\limits_{\bfO(h)}g(t) dt,$$
and 
so we have an Euler Maclaurin formula for regular orthants:
\begin{equation}
\labell{EMregularorthant}
{\sum_{{\bf C}\cap \Z^n}}'f=
\left. \prod_{i=1}^n\Ltk(\deldel{h_i})
\int\limits_{{\bf C}(\pm h_1,\dots,\pm h_n)}f(x)dx
\right|_{h_1=\cdots=h_n=0} +R^{\bf C}_{m}(f),
\end{equation}
where
$$R^{\bf C}_{m}(f)=R^{\st}_{m}(g).$$

\bigskip

Let $\Delta$ be a compact convex polytope.
As in Section \ref{sec:decompose},
$\Delta$ can be written as an intersection of half-spaces
\begin{equation} \labell{Delta}
\Delta = H_1 \cap \ldots \cap H_d, \quad \text{where} \quad
H_i = \{ x \mid \langle u_i , x \rangle + \mu_i \geq 0 \}
\quad \text{ for } i=1,\ldots,d
\end{equation}
and $d$ is the number of facets of $\Delta$.
The vector $u_i \in {\R^n}^*$ can be thought of as the inward normal
to the $i$th facet of $\Delta$;
a-priori it is determined up to multiplication by a positive number.
If all the vertices of $\Delta$ are integral, then the $u_i$'s can
be chosen to belong to the dual lattice ${\Z^n}^*$, and we can 
fix our choice of the $u_i$'s by imposing
the normalization condition that the $u_i$'s be primitive lattice elements,
that is, that no $u_i$ can be expressed as a multiple of a lattice
element by an integer greater than one.
(The fact that a normal vector $u$ to a facet $\sigma$
can be chosen to be integral is a consequence of Cramer's rule.
Indeed, we can choose integral edge vectors
$\beta_1,\ldots,\beta_n$ that emanate from a vertex on $\sigma$ 
such that $\beta_1,\ldots\beta_{n-1}$ span the tangent plane to $\sigma$ 
and $\beta_n$ is transverse to $\sigma$. 
Solving the linear equations $\left< u , \beta_1 \right> =$
$\ldots$ $=\left< u , \beta_{n-1} \right> = 0$ and
$\left< u , \beta_{n} \right> = 1$,
we get an inward normal vector $u$ with rational entries;
clearing denominators, we may assume that $u$ is actually integral.)

We can then consider the ``dilated polytope" $\Delta(h_1,\ldots,h_d)$,
which is obtained by shifting the $i$th facet outward
by a ``distance" $h_i$. More precisely,
$$ \Delta(h) = \bigcap_{i=1}^d \{ x \mid \langle u_i , x \rangle 
   + \mu_i + h_i \geq 0 \} 
\quad \text{where} \quad h = (h_1,\ldots,h_d). $$

\bigskip

Now assume that $\Delta$ is simple.
Then $\Delta(h)$ is simple if $h$ is sufficiently small.
The polar decomposition of $\Delta(h)$ involves ``dilated orthants".
However, dilating the facets of $\Delta$ outward results in dilating
some facets of $\bfC_v^\sharp$ inward and some outward.  
Explicitly, for $i \in I_v = \{ i_1, \ldots, i_n \} $, 
the inward normal vector to the $i$th facet of $\bfC_v^\sharp$ is 
\begin{equation} \labell{u i v sharp}
   u_{i,v}^\sharp = \begin{cases} 
 u_i & \text{ if } \alpha_{i,v}^\sharp =  \alpha_{i,v} \\
-u_i & \text{ if } \alpha_{i,v}^\sharp =  -\alpha_{i,v} .
\end{cases}
\end{equation}
Hence, the dilated orthants that occur on the right hand side
of the polar decomposition of $\Delta(h)$
are $\bfC_v^\sharp(h_{i_1,v}^\sharp, \ldots, h_{i_n,v}^\sharp)$,
where
\begin{equation} \labell{h i v sharp}
 h_{i,v}^\sharp = \begin{cases}
 h_i & \text{ if } \alpha_{i,v}^\sharp =  \alpha_{i,v}  \\
-h_i & \text{ if } \alpha_{i,v}^\sharp = - \alpha_{i,v}. \end{cases}
\end{equation}
This subtlety in the signs does not effect the formula for regular polytopes,
because of the symmetry of $\Ltk$; however, it is needed to derive
the formula for simple polytopes.

If $\Delta$ is a regular integral polytope,
then the $C_v^\sharp$'s are regular integral orthants,
and the ``dilated orthants" are exactly as in \eqref{uiss}.
We then have, by \eqref{wheee}, \eqref{EMregularorthant},
and the symmetry of $\Ltk$,
\begin{multline}
\labell{cow}
{\sum_{\Delta\cap \Z^n}}'f
 = \sum_v (-1)^{\#v}{\sum_{\bfC_v^\sharp \cap \Z^n}}'f \\
 = \sum_v(-1)^{\#v}\left(\left. \prod_{i\in I_v=\{i_1,\dots,i_n\}}
\Ltk(\deldel{h_i})
\int\limits_{\bfC_v^\sharp (\pm h_{i_1},\dots,\pm h_{i_n})}f(x)dx
\right|_{h=0} +R^{\bfC_v^\sharp}_{m}(f)\right).
\end{multline}
We may multiply the differential operator $\prod \Ltk(\deldel{h_i})$
in the above expression by any number of operators of the form
$\Ltk(\deldel{h_j})$ where $j \not\in I_v$,
since all that will remain of these operators is
the constant term $1$, all actual differentiations yielding zero.
The right hand side of \eqref{cow} is then equal to
\begin{multline*}
 \left. \prod_{i=1}^d \Ltk(\deldel{h_i}) \sum_v (-1)^{\# v} 
\int\limits_{\bfC_v^\sharp (h_{i_1,v}^\sharp,\dots,h_{i_n,v}^\sharp)} f(x)dx 
  \right|_{h=0} + R^{\bfC_v^\sharp}_{m}(f) \\
 = \left. \prod_{i=1}^d\Ltk(\deldel{h_i})\int\limits_{{\Delta}(h_1,\dots,h_d)}
f(x)dx\right|_{h=0}+
S_\Delta^{m}(f)
\end{multline*}
where
\begin{equation}
\labell{polytoperemainder}
S_\Delta^{m}(f):=\sum_v(-1)^{\#v}R^{\bfC_v^\sharp}_{m}(f)
\end{equation}
and where $\{ i_1, \ldots, i_n \} = I_v$ in the $v$'th summand.

Notice that both
$${{\sum_{\Delta\cap \Z^n}}'f}$$
and
$$\left. \prod_{i=1}^d\Ltk(\deldel{h_i})\int\limits_{{\Delta}(h_1,\dots,h_d)}
f(x)dx\right|_{h=0}$$
do not
depend on the choice of polarization, and both vanish on any function $f$
whose support
is disjoint from the polytope. So the same must be true of the remainder.
We have proved the following result:

\begin{Theorem}[\cite{PNAS}] \labell{EMregular} 
Let $m \geq 1$ be an integer.
Let $\Delta \subset \R^n$ be a regular integral
polytope and $f$ a $C^{mn}$ function
of compact support on $\R^n$. Choose a ``polarizing vector" for $\Delta$.
Then
$${\sum_{\Delta\cap \Z^n}}'f=
\left. \prod_{i=1}^d\Ltk(\deldel{h_i})\int\limits_{{\Delta}(h_1,\dots,h_d)}
f(x)dx\right|_{h=0}+ S_\Delta^{m}(f)$$
where $k = \lfloor m/2 \rfloor$ and
where $S_\Delta^{m}(f)$ is given by \eqref{polytoperemainder}.
This remainder can be expressed as a sum of integrals over orthants 
of bounded periodic functions times various partial derivatives
of $f$ of order no less than $m$ and no more than $mn$.  Finally,
this remainder is independent of the choice of
polarization and is a distribution supported on the polytope $\Delta$.
\end{Theorem}

% -----------------------------------------------------------------------------
\section{Finite groups associated to a simple integral polytope and its faces.}
% -----------------------------------------------------------------------------
\labell{sec:groups}

In order to extend Theorem \ref{EMregular} to simple integral polytopes
that may not be regular, we must extend the Euler Maclaurin
formula for regular orthants (\ref{EMregularorthant}) to a
formula that is valid for simple orthants that may not be regular.
In this section we analyze
certain finite groups that arise in this generalization.

Let $\bfC$ be a simple integral orthant. 
This means that we can write $\bfC$ as the intersection of $n$ half-planes
in general position,
\begin{equation} \labell{general position}
 \bfC = H_1 \cap \ldots \cap H_n \quad \text{where} \quad
  H_i = \{ x \mid \langle u_i , x \rangle + \mu_i \geq 0 \}
\quad \text{ for } i=1,\ldots,n,
\end{equation}
and that the vertex of $\bfC$ is in $\Z^n$.
The $u_i$'s are inward normals to the facets of $\bfC$,
and we choose them to be primitive elements
of the dual lattice ${\Z^n}^*$. (See Section \ref{sec:EM regular}.)
We choose $\alpha_1,\ldots,\alpha_n$ to be the dual basis,
that is,
$$ \l< u_j , \alpha_i \r> = \begin{cases}
1 & j=i \\ 0 & j \neq i.
\end{cases}$$
The $\alpha_i$'s are edge vectors of $\bfC$, but they might not be integral:
they generate a lattice in $\R^n$ which is a finite extension of $\Z^n$.
This extension is trivial exactly if $\Delta$ is regular at $v$,
that is, if the $u_i$'s, generate the dual lattice ${\Z^n}^*$. 
To the cone $\bfC$ we associate the finite group
\begin{equation} \labell{Gamma}
 \Gamma := {\Z^n}^* / \sum \Z u_i .
\end{equation}
So this group is trivial exactly if the cone $\bfC$ is regular.

\begin{Lemma} \labell{character}
In the above setting,
\begin{equation} \labell{the character}
 \gamma \mapsto e^{2 \pi i \left< \gamma , x \right> } 
\end{equation}
is a well defined character on $\Gamma$ whenever 
$x = m_1 \alpha_1 + \ldots + m_n \alpha_n$ where $m_j$ are integers.
This character is trivial if and only if $x \in \Z^n$.
\end{Lemma}

\begin{proof}
Let $ \tilde{\gamma} \in {\Z^n}^*$ be an element that represents $\gamma$.
If we expand it as a combination of the basis elements $u_j$,
so that
$$ \tilde{\gamma} = b_1 u_1 + \ldots + b_n u_n,$$
then $\left< \tilde{\gamma} , \alpha_j \right> = b_j$ for all $j$.
If $\tilde{\gamma}'$ is another element of ${\Z^n}^*$ 
that represents $\gamma$, then it differs from $\tilde{\gamma}$ 
by an integral combination of the $u_i$'s. So
$\left< \tilde{\gamma}' , \alpha_j \right> $
differs from $ \left< \tilde{\gamma} , \alpha_j \right> $
by an integer.
Hence, when $x$ is an integral combination of the $\alpha_j$'s,
the pairing
$\left< \tilde{\gamma} , x \right> $ is well defined modulo $\Z$,
and so $e^{2 \pi i \left< \gamma,x \right>}$ is well defined.

Finally, the character \eqref{the character} is trivial
if and only if $\left< \tilde{\gamma} , x \right> $ is an integer
for all $\gamma \in {\Z^n}^*$, and this holds if and only if $x \in \Z^n$.
\end{proof}

\bigskip

Let $\Delta$ be a simple integral polytope in $\R^n$,
given by \eqref{Delta}.
For each vertex $v$ of $\Delta,$ let $I_v \subset \{ 1 , \ldots , d \}$
denote the set of facets of $\Delta$ that meet at $v$.
The normal vectors $u_i$,
$i \in I_v$, form a basis of ${\R^n}^*$.
We choose
\begin{equation} \labell{alpha i v}
\valpha_{i,v} \in \R^n \quad , \quad i \in I_v
\end{equation}
to be the dual basis.

Given any face $F$ of the polytope $\Delta$,
let $I_F$ denote the set of facets of $\Delta$ which meet at $F$.
Because $\Delta$ is simple, the vectors $u_i$, for $i \in I_F$,
are linearly independent.
Let $N_F \subseteq {\R^n}^*$ be the subspace
$$ N_F = \span \{ u_i \mid i \in I_F \}.$$

\begin{Remark}
It is natural to define the tangent space to the face $F$ to be
$TF = \span \{x-y \mid x,y \in F\}$,
the normal space to $F$ to be the quotient $\R^n / TF$,
and the co-normal space to be the dual of the normal space.
With these definitions, $N_F$ is the co-normal space to $F$.
\end{Remark}

To each face $F$ of $\Delta$ we associate a finite abelian
group $\Gamma_F$. Explicitly, the lattice
$$\vV_F = \sum_{i \in I_F} \Z u_i \, \subset \, N_F $$
is a sublattice of $N_F \cap {\Z^n}^*$ of finite index, and
the finite abelian group associated to the face $F$ is
the quotient
\begin{equation} \labell{def GammaF}
\Gamma_F := (N_F \cap {\Z^n}^*) / \vV_F.
\end{equation}
If $F=v$ is a vertex, this is the same as the finite abelian group
associated to the tangent cone $\bfC_v$ as in \eqref{Gamma}.

Let $E$ and $F$ be two faces of $\Delta$ with $F \subseteq E$.
This inclusion implies that $I_E \subseteq I_F$, and hence
\begin{equation} \labell{uis}
 \{ u_i \}_{i \in I_E} \subseteq \{ u_i \}_{i \in I_F}.
\end{equation}
Because these sets are bases of the vector spaces $N_E$ and $N_F$,
we have an inclusion
$$N_E \subseteq N_F.$$
Because the sets occurring in \eqref{uis} 
are $\Z$-bases of the lattices $\vV_E$ and $\vV_F$,
we have $N_E \cap \vV_F = \vV_E$.
Hence, the natural map from  $\Gamma_E = ({\Z^n}^* \cap N_E) / \vV_E$
to $\Gamma_F = ({\Z^n}^* \cap N_F) / \vV_F$ is one to one,
and it provides us with a natural inclusion map:
$$ \text{if \  $F \subseteq E$ \  then \  $\Gamma_E \subseteq \Gamma_F$.} $$

We define a subset $\Gamma_F^\flat$ of $\Gamma_F$ by 
\begin{equation} \labell{def:GammaFsharp}
\Gamma_F^\flat
  := \Gamma_F \ssminus
\bigcup_{\text{faces } E \text{ such that }  E \supsetneq F} \Gamma_E.
\end{equation}
Then
\begin{equation} \labell{Gammas}
 \Gamma_v = \bigsqcup_{\{F:v \in F\}}  \Gamma_F^\flat.
\end{equation}

Recall that 
\begin{equation} \labell{lambda j v}
 \lambda_{\gamma,j,v} := 
   e^{ 2 \pi i \left< \gamma ,\alpha_{j,v} \right> } , 
   \quad \text{ for } \gamma \in \Gamma_v
   \quad \text{ and } j \in I_v,
\end{equation}
is well defined, by Lemma \ref{character}.  It is a root of unity.
We will need the following results.

\begin{Claim} \labell{claim1}
If $\gamma \in \Gamma_F$ and $j \in I_F$,
then $\lambda_{\gamma,j,v}$ is the same for all $v \in F$.
\end{Claim}

This allows us to define $\lambda_{\gamma,j,F}$
for $\gamma \in \Gamma_F$ and $j \in I_F$ such that
$$ \lambda_{\gamma,j,F} = \lambda_{\gamma,j,v}
\quad \text{ for } \gamma \in \Gamma_F  
\text{ and } j \in I_F, \text{ if } v \in F .$$

\begin{Claim} \labell{claim2}
If $\gamma \in \Gamma_F$ and $j \in I_v \ssminus I_F$
then $\lambda_{\gamma,j,v}$ is equal to one.
\end{Claim}

This allows us to define
$\lambda_{\gamma,j,F} = 1$ when $\gamma \in \Gamma_F$
and $j \in \{ 1 , \ldots, d \} \ssminus I_F $.
Then
\begin{equation} \labell{horse}
 \lambda_{\gamma,j,F} = \lambda_{\gamma,j,v}
\quad \text{ for } \gamma \in \Gamma_F  
\text{ and } 1 \leq j \leq d , \text{ if } v \in F 
\end{equation}
and
\begin{equation} \labell{lambda is one}
\lambda_{\gamma,j,F} = 1 \quad \text{ for } \gamma \in \Gamma_F
\text{ if } j \not\in I_F.
\end{equation}

\begin{Claim} \labell{claim3}
If $\gamma \in \Gamma_F^\flat$ and $j \in I_F$, 
then $ \lambda_{\gamma,j,F} \neq 1$.
\end{Claim}

\begin{proof}[Proof of Claims \ref{claim1}, \ref{claim2}, and \ref{claim3}]
Let $\gamma \in \Gamma_F$ be represented by
$$ \tilde{\gamma} = 
   \sum_{i \in I_F} b_i u_i \in N_F \cap {\Z^n}^* $$
for some $b_i \in \R$.  (See \eqref{def GammaF}.)

Let $v \in F$.  Because $\alpha_{j,v}$, $j \in I_v$, 
is a dual basis to $u_j$, $j \in I_v$, we have
$$ \left< \tilde{\gamma} , \alpha_{j,v} \right>
 = \begin{cases}  b_j & j \in I_F \\ 0 & j \in I_v \ssminus I_F.
\end{cases}$$
Hence,
$$ \lambda_{\gamma,j,v} 
   = e^{2 \pi i \left< \tilde{\gamma} , \alpha_{j,v} \right> }
   = \begin{cases} e^{ 2 \pi i b_j} & j \in I_F \\ 
     1 & j \in I_v \ssminus I_F \end{cases}
$$
is independent of $v$ and is equal to $1$ 
if $j \in I_v \ssminus I_F$.
This prove Claims \ref{claim1} and \ref{claim2}.

Let $j \in I_F$. If
$ \lambda_{\gamma,j,F} := e^{ 2 \pi i b_j } $
is equal to one, then $b_j$ is an integer, so
\begin{equation} \labell{in GammaE}
 \tilde{\gamma}' = 
 \sum_{i \in I_F \ssminus \{ j \} } b_i u_i 
\end{equation}
also represents $\gamma$.  Let $E \supset F$
be the face such that $I_E = I_F \ssminus \{ j \}$.
Then, by \eqref{in GammaE}, $\gamma \in \Gamma_E$.
This proves Claim \ref{claim3}.
\end{proof}

% -----------------------------------------------------------------------------
\section{Euler Maclaurin formula with remainder for simple integral polytopes.}
% -----------------------------------------------------------------------------
\labell{sec:EM simple}

Let us begin by deriving an Euler Maclaurin formula with remainder
for a simple integral orthant. 

We recall the set-up of Section \ref{sec:groups}.
Let $\bfC$ be a simple integral orthant.  Let $u_1, \ldots, u_n$
be the inward normals to its facets, chosen to be primitive
elements of the dual lattice ${\Z^n}^*$, 
and let $\alpha_1 , \ldots, \alpha_n$
be the dual basis to the $u_i$'s, so that 
$$ \bfC = v + \sum_{j=1}^n \R_{\nonneg} \alpha_j.$$
Let
$$ \Gamma = {\Z^n}^* / \sum \Z u_j $$
be the finite group associated to $\bfC$.
By Lemma \ref{character},
$\gamma \mapsto e^{ 2 \pi i \left< \gamma, x \right> }$
defines a character on $\Gamma$
whenever $x \in \sum \Z \alpha_j$,
and this character is trivial if and only if $x \in \Z^n$.
By a theorem of Frobenius, the average value of a character
on a finite group is zero if the character is non-trivial
and one if the character is trivial.  So
$$ \frac{1}{| \Gamma |} \sum_{\gamma \in \Gamma}
e^{ 2 \pi i \left< \gamma, x \right> } = \begin{cases}
 1 & \text{if } x \in \Z^n \\
 0 & \text{if } x \not \in \Z^n \\
\end{cases} $$
for all $x \in \sum \Z \alpha_j$.
Then, for a compactly supported function $f(x)$ on $\R^n$,
\begin{eqnarray}
\nonumber
 {\sum_{\bfC \cap \Z^n}}{'} f
 & = & {\sum_x}' \left( \frac{1}{|\Gamma|}
\sum_{\gamma \in \Gamma} e^{2\pi i \left< \gamma , x \right>} \right) f(x) \\
\labell{eqAA}
 & = & \frac{1}{|\Gamma|} \sum_{\gamma \in \Gamma} {\sum_x}{'} 
       e^{2\pi i \left< \gamma , x \right>} f(x) 
\end{eqnarray}
where we sum over all 
\begin{equation} \labell{x}
 x = v + m_1 \alpha_1 + \ldots + m_n \alpha_n ,
\end{equation}
where the $m_i$'s are non-negative integers.

The cone $\bfC$ is the image of the standard orthant $\bfO$
under the affine map
\begin{equation} \labell{affine map}
 (t_1,\ldots,t_n) \mapsto x = v + t_1 \alpha_1 + \ldots + t_n \alpha_n .
\end{equation}
This map sends the lattice $\Z^n$ onto the lattice
$\sum \Z \alpha_j$.    The inverse transformation is given by
\begin{equation} \labell{inverse transformation}
 t_i = \left< u_i , x-v \right>  .
\end{equation}

Let us concentrate on one element $\gamma \in \Gamma$.
Because $v \in \Z^n$, from \eqref{x} we get
$$ e^{2\pi i \left< \gamma , x \right>} 
   = \prod_{j=1}^n \lambda_{j} ^{m_j} 
\quad \text{where} \quad
\lambda_{j} = e^{2\pi i\left< \gamma , \alpha_j \right> } ,$$
so that the inner sum in (\ref{eqAA}) becomes
\begin{equation} \labell{eqBB}
{\sum_{m_1 \geq 0}}' \lambda_{1}^{m_1} \cdots
{\sum_{m_n \geq 0}}' \lambda_{n}^{m_n} g(m_1,\ldots,m_n),
\end{equation}
where
$$ g(t_1,\ldots,t_n) = 
   f(v + t_1 \alpha_1 + \ldots + t_n \alpha_n).  $$

Recall that we had the twisted remainder formula
$$ {\sum_{m \geq 0}}' \lambda^m g(m) 
 = \left. \bfM^{k,\lambda}(\deldel{h}) \int_{-h}^\infty g(t) dt \right|_{h=0}
 + (-1)^{k-1} \int_0^\infty Q_{k,\lambda}(t) g^{(k)}(t) dt$$
for all compactly supported functions $g(x)$ of type $C^k$,
where $k \geq 1$, where $\lambda$ is a root of unity,
and where $M^{k,\lambda}$ is a polynomial of degree $\leq k$.
(See \eqref{twistedemrayeq} and Remark \ref{M vs L}.)

Iterating this formula, the sum in \eqref{eqBB} can be written as
$$
  \bfM^{k,\lambda_{1}}(\deldel{h_1}) \int_{-h_1}^\infty \cdots
     \bfM^{k,\lambda_{n}}(\deldel{h_n}) \int_{-h_n}^\infty 
     g(t_1,\ldots,t_n) dt_1 \cdots dt_n 
    + R^{st}_k(\lambda_{1},\dots,\lambda_{n};g)
$$
\begin{equation} \labell{integral over Oh}
     = \prod_{i=1}^n \bfM^{k,\lambda_i}(\deldel{h_i})
       \int\limits_{\bfO(h)} g(t_1,\ldots,t_n) dt_1 \cdots dt_n 
      + R^{st}_k(\lambda_{1},\ldots,\lambda_{n};g)
\end{equation}
with
$$ \bfO(h) = \{ (t_1,\ldots,t_n) \ | \ t_i \geq -h_i \text{\ for all $i$} \}$$
and where the remainder is given by
\begin{multline*}
R^{st}_k(\lambda_1,\dots,\lambda_n ; g) = \\
\sum_{I \subsetneq \{ 1, \ldots, n \} } (-1)^{(k-1)(n-|I|)} 
\prod_{i \in I} \bfM^{k,\lambda_i}(\deldel{h_i}) 
\int\limits_{O(h)} 
\left.
\prod_{i \notin I} Q_{k,\lambda_i} (t_i)
\prod_{i \notin I} \frac{\del^{k}}{\del t_i ^{k}}
 g(t_1,\ldots,t_n) dt_1 \cdots dt_n
\right|_{h=0} .
\end{multline*}
Using \eqref{formula with varphi} to
express this as a sum of integrals over the (non-shifted) orthant,
we get
\begin{multline} \labell{remainder}
R^{st}_k(\lambda_1,\dots,\lambda_n ; g) = \\
\sum_{I \subsetneq \{ 1, \ldots, n \} } (-1)^{(k-1)(n-|I|)} 
\int\limits_{\bfO}
\prod_{i \in I} \bfM^{k,\lambda_i}( - \deldel{t_i}) 
\prod_{i \notin I} Q_{k,\lambda_i} (t_i)
\prod_{i \notin I} \frac{\del^{k}}{\del t_i ^{k}}
 g(t_1,\ldots,t_n) dt_1 \cdots dt_n .
\end{multline}

We now perform in \eqref{integral over Oh} and \eqref{remainder}
the change of variable given by the transformation
\eqref{inverse transformation}.
The integrals over $\bfO(h)$ and $\bfO$
get replaced by integrals over $\bfC(h)$ and $\bfC$
times the Jacobian factor $|\Gamma|$,
the function $g(t)$ gets replaced by the function $f(x)$,
and the partial derivative $\deldel{t_i}$ gets replaced by
the directional derivative $D_{\alpha_i}$.
Summing the expressions \eqref{integral over Oh} over $\gamma \in \Gamma$
and dividing by $|\Gamma|$, we get, by \eqref{eqAA},
\begin{equation} \labell{cat}
 \sum_{\bfC \cap \Z^n}{'} f = \sum_{\gamma \in \Gamma}
\prod_{i=1}^n \bfM^{k,\lambda_{\gamma,i}}(\deldel{h_i}) 
\left. \int_{\bfC(h)} f(x) dx \right|_{h=0} + R_k^{\bfC} (f) 
\end{equation}
where
$$ \lambda_{\gamma,j} := e^{2 \pi i \left< \gamma, \alpha_j \right> }$$
and where
\begin{multline} \labell{remainder C}
 R^{\bfC}_k(f) := \\
\sum_{\gamma \in \Gamma}
     \sum_{I \subsetneq \{ 1, \ldots, n \} } (-1)^{(k-1)(n-|I|)}
     \int\limits_{\bfC}
     \prod_{i \in I} \bfM^{k,\lambda_{\gamma,i}} (-D_{\alpha_i})
     \prod_{i \not \in I}
               Q_{k,\lambda_{\gamma,i}} ( \left< u_i , x-v \right> )
     \prod_{i \not \in I} (D_{\alpha_i})^k f(x) dx .
\end{multline}

\bigskip

Let $\Delta$ be a simple polytope, given by \eqref{Delta}.
Choose a polarizing vector for $\Delta$ and let $\bfC_v^\sharp$
denote the polarized tangent cones.
The inward normals to the facets of $\bfC_v^\sharp$ are 
given by \eqref{u i v sharp}, the dual basis to these vectors 
is $\alpha_{i,v}^\sharp$, $i \in I_v$,
and the roots of unity that appear in the Euler Maclaurin formula
for $\bfC_v^\sharp$ are then
\begin{equation} \labell{lambda sharp}
\lambda_{\gamma,i,v}^\sharp 
= e^{2 \pi i \langle \gamma, \alpha_{i,v}^\sharp \rangle } = \begin{cases} 
 \lambda_{\gamma,i,v} & \text{ if } \alpha_{i,v}^\sharp =  \alpha_{i,v} \\
\lambda_{\gamma,i,v}\inv & \text{ if } \alpha_{i,v}^\sharp =  -\alpha_{i,v} .
\end{cases}
\end{equation}
Also recall that the polar decomposition of $\Delta(h)$ 
involves the dilated orthants 
$C_v^\sharp(h_{i_1,v}^\sharp, \ldots, h_{i_n,v}^\sharp )$
where $h_{i,v}^\sharp$ are as in \eqref{h i v sharp}.

Let $k \geq 1$ be an integer.
For any compactly supported function $f$ on $\R^n$ of type $C^{nk}$,
the polar decomposition of $\Delta(h)$ and the formula \eqref{cat} 
give
\begin{multline} \labell{sum for simple}
{\sum_{\Delta \cap \Z^n}}' f = 
\sum_v (-1)^{\# v} {\sum_{\bfC_v^\sharp \cap \Z^n}}' f \\
 = \sum_v (-1)^{\# v} \sum_{\gamma \in \Gamma_v}
   \prod_{i \in I_v = \{ i_1, \ldots, i_n \} } 
   \bfM^{k,\lambda_{\gamma,i,v}^\sharp}(\deldel{h_{i,v}^\sharp}) \left.
   \int_{\bfC_v^\sharp(h_{i_1,v}^\sharp,\ldots,h_{i_n,v}^\sharp)}
   f(x) dx \right|_{h=0} 
   + R_k^{\Delta}(f),
\end{multline}
where the remainder is given by 
\begin{equation} \labell{remdef}
 R^\Delta_k(f) := \sum_v (-1)^{\#v} R^{C_v^\sharp}_k(f) .
\end{equation}
Note that either $h_{i,v}^\sharp = h_{i}$  and 
$\lambda_{\gamma,i,v}^\sharp = \lambda_{\gamma,i,v}$,
or $h_{i,v}^\sharp = -h_{i}$  and 
$\lambda_{\gamma,i,v}^\sharp = \lambda_{\gamma,i,v}\inv$.
By the symmetry property \eqref{Dandlambdainverse},
this gives
$$ \bfM^{k,\lambda_{\gamma,i,v}^\sharp} (\deldel{h_{i,v}^\sharp})
 = \bfM^{k,\lambda_{\gamma,i,v}} (\deldel{h_i}) .$$

For $j \not\in I_v$, because $\lambda_{\gamma,j,v} = 1$
(see \eqref{lambda is one}),
we have
$\bfM^{k,\lambda_{\gamma,j,v} }(\deldel{h_j}) = 1 + $powers 
of $\deldel{h_j}$.
Also, the cone 
$C_v^\sharp( h_{i_1,v}^\sharp, \ldots, h_{i_n,v}^\sharp)$ 
is independent of $h_j$ for $j \not \in I_v$.
Therefore, \eqref{sum for simple} is further equal to
\begin{equation} \labell{eqA}
\left. \sum_v (-1)^{\# v} \sum_{\gamma \in \Gamma_v} 
   \prod_{j=1}^d \bfM^{k,\lambda_{\gamma,j,v}} (\deldel{h_j}) 
   \int\limits_{\bfC_v^\sharp(h_{i_1,v}^\sharp, \ldots, h_{i_n,v}^\sharp)} 
   f(x) dx \right|_{h=0} + R^\Delta_k(f)
\end{equation}
where in the $v$'th summand $\{ i_1, \ldots, i_n \} = I_v$. 

Define
\begin{equation}
\labell{def of Mk gamma F}
\bfM^k_{\gamma,F} = \prod_{j=1}^d \bfM^{k,\lambda_{\gamma,j,F}}(\deldel{h_j})
\quad \text{ for } \gamma \in \Gamma_F .
\end{equation}
Then we have, by \eqref{horse},
\begin{equation} \labell{Fv}
\bfM^k_{\gamma,F} = \bfM^k_{\gamma,v} \quad \text{ whenever } 
\gamma \in \Gamma_F \text{ and } v \in F,
\end{equation}
where we identify $\gamma \in \Gamma_F$ with its image
under the inclusion map $\Gamma_F \hookrightarrow \Gamma_v$.

Then \eqref{eqA} is equal to
\begin{eqnarray}
\nonumber
\left.\sum_v (-1)^{\# v} \sum_{\gamma \in \Gamma_v}
\bfM_{\gamma,v}^k \int_{C_v^\sharp(h_{i_1,v}^\sharp,\ldots,h_{i_n,v}^\sharp)} 
    f(x) dx \right|_{h=0} + R^\Delta_k(f) \\
\labell{eqB}
  =  \left.\sum_F \sum_{\gamma \in \Gamma_F^\flat}
       \bfM_{\gamma,F}^k \sum_{v \in F} (-1)^{\# v}
       \int_{C_v^\sharp(h_{i_1,v}^\sharp,\ldots,h_{i_n,v}^\sharp)} 
       f(x) dx \right|_{h=0} + R^\Delta_k(f) ,
\end{eqnarray}
by \eqref{Gammas} and \eqref{Fv}.
In the interior summation we may now add similar summands that correspond
to $v \not \in F$.  
These summands make a zero contribution to (\ref{eqB})
for the following
reason.  If $v \not \in F$ then there exists $i \in I_F \ssminus I_v$.
Because $i \not \in I_v = \{ i_1, \ldots, i_n \}$, the cone
$C_v^\sharp(h_{i_1,v}^\sharp , \ldots h_{i_n,v}^\sharp)$ 
is independent of $h_i$.
So it is enough to show that $\bfM_{\gamma,F}^k$ is a multiple 
of $\deldel{h_i}$.
But because $\gamma \in \Gamma_F^\flat$ and $i \in I_F$, we have
$\lambda_{\gamma,i,F} \neq 1$. (See Claim \ref{claim3}.)
By Remark \ref{M vs L}, this
implies that $\bfM^{k,\lambda_{\gamma,i,F}}(\deldel{h_i})$,
which is one of the factors in $\bfM_{\gamma,F}^k$, is a multiple 
of $\deldel{h_i}$.  Hence, \eqref{eqB} is equal to
\begin{multline}
\left.\sum_F \sum_{\gamma \in \Gamma_F^\flat}
       \bfM_{\gamma,F}^k \sum_{\text{all } v} (-1)^{\# v}
       \int_{C_v^\sharp( h_{i_1,v}^\sharp ,\ldots, h_{i_n,v}^\sharp )} 
       f(x) dx \right|_{h=0} + R^\Delta_k(f)  \\
     = \left.\sum_F \sum_{\gamma \in \Gamma_F^\flat}
       \bfM_{\gamma,F}^k \int_{\Delta(h)} f(x) dx \right|_{h=0}
       + R^\Delta_k(f) .\\
\end{multline}

We have therefore proved our main result:

\begin{Theorem}\labell{main}
Let $\Delta$ be a simple integral polytope in $\R^n$ and let
$f\in C^{nk}_c(\R^n)$ be a compactly supported function on $\R^n$
for $k \geq 1$.  Choose a polarizing vector for $\Delta$.  Then
$${\sum_{\Delta \cap \Z^n}}' f = \left.\sum_F \sum_{\gamma \in \Gamma_F^\flat}
       \bfM_{\gamma,F}^k \int_{\Delta(h)} f(x) dx \right|_{h=0}
       + R^\Delta_k(f) $$
where $\bfM_{\gamma,F}^k$ are differential operators
defined in \eqref{def of Mk gamma F} and where the
remainder $R^\Delta_k(f)$ is given by equation (\ref{remdef}).
Moreover, the differential operators $\bfM_{\gamma,F}^k$
are of order $\leq k$ in each of the variables $h_1,\ldots,h_d$.
Also, the remainder is a sum of integrals
over orthants of bounded periodic functions times various partial derivatives 
of $f$ of order no less than $k$ and no more than $kn$.
Finally, this remainder is independent of the choice of polarization
and is a distribution supported on the polytope $\Delta$.
\end{Theorem}

% =========================================================
\section{Estimates on the remainder and 
         Euler Maclaurin formulas for symbols and
	 for polynomials}
\labell{sec:estimates}
% =========================================================

We first recall a definition from the theory of partial
differential equations.  A smooth function $f \in C^\infty(\R^n)$ 
is called a \textbf{symbol of order} $\mathbf{N}$
if for every $n$-tuple of non-negative integers
$a :=(a_1,\dots,a_n)$,
there exists a constant $C_a$ such that 
$$|\partial_1^{a_1}\dots \partial_n^{a_n} f(x) | \leq C_a (1 + |x|)^{N - |a|}$$
where $|a| = \sum_i a_i$.  In particular, a polynomial of degree
$N$ is a symbol of order $N$.  
Note that if $f$ is a symbol of order $N$ on $\R^n$ then its
derivatives of order $a$ are in $L^1$ if $N < |a| - n$.
In this section we will
show that the Euler Maclaurin formula of Theorem \ref{main}
can be extended to symbols, and gives rise in this way
to an {\em exact} Euler Maclaurin formula for polynomials.

To make further progress,
we first require an estimate on the remainder term $R^\Delta_k(f)$.
We recall (see Theorem \ref{main}) that this remainder
can be expressed as a sum of integrals over orthants
of bounded periodic functions time various partial derivatives of $f$
of order no less than $k$ and no more than $kn$.
Explicitly, from \eqref{remdef} and \eqref{remainder C} we get
\begin{multline} \labell{def of RkDelta}
 R_k^\Delta (f) = \sum_v (-1)^{\# v} \sum_{\gamma \in \Gamma_v}
     \sum_{I \subsetneq I_v} (-1)^{(k-1)(n-|I|)} \\
     \int\limits_{\bfC_v^\#}
     \prod_{j \in I} 
         \bfM^{k,\lambda_{\gamma,j,v}^\sharp} (-D_{\alpha_{j,v}^\sharp})
     \prod_{j \not \in I}
        Q_{k,\lambda_{\gamma,j,v}^\sharp} 
                         (\langle u_{j,v}^\sharp, x-v \rangle )
     \prod_{j \not \in I} (D_{\alpha_{j,v}^\sharp})^k
     f(x) dx . 
\end{multline}

The functions $Q_{k,\lambda}$ are bounded and periodic,
and $M^{k,\lambda}(-D_{\alpha_j})$ are differential operators
of order $k$.  In particular, each integrand in the formula
for $R_k^\Delta(f)$ is dominated by a constant times
$$ \sup_{ \{ j_1,\ldots,j_n \} } |\del_1^{j_1} \cdots \del_n^{j_n} f| $$
where the supremum is taken over all $n$-tuples $\{ j_1,\ldots,j_n \}$
with $k \leq j_1 + \ldots + j_n \leq nk$. 
Consequently, $R_k^\Delta(f)$ is well defined (by the same formula
\eqref{def of RkDelta}) for any smooth function $f$ whose derivatives 
of order between $k$ and $nk$ are integrable on $\R^n$.  
In particular, it is defined when $f$ is a symbol of order less than $k-n$.
Moreover, we get an estimate for the remainder:

\begin{Proposition} \labell{polyest}
The distribution $R_k^\Delta(\cdot)$ extends to symbols of order
less than $k-n$, by the explicit expression \eqref{def of RkDelta} above,
and satisfies an estimate of the form
$$|R^\Delta_k(f)| 
\leq K(k,\Delta) \cdot {\rm sup}_{\{j_1,\dots,j_n\}} 
|\partial_1^{j_1}\dots\partial_n^{j_n} f|_{L_1(\R^n)},$$
where the supremum is taken over all $n$-tuples $\{j_1,\cdots,j_n\}$
with $k \leq j_1 + \cdots + j_n  \leq nk.$
\end{Proposition}

By applying this estimate we will obtain the following Euler Maclaurin formula
for symbols.  

\begin{Theorem}\labell{syms}
Let $\Delta$ be a simple integral polytope in $\R^n$,
let $f$ be a symbol of order $N$ on $\R^n$, and choose $k \geq N + n + 1$.
Then 
\begin{equation} \labell{EM formula for symbols}
{\sum_{\Delta \cap \Z^n}}' f = 
\left.\sum_F \sum_{\gamma \in \Gamma_F^\flat}
       \bfM_{\gamma,F}^k \int_{\Delta(h)} f(x) dx \right|_{h=0}
       + R^\Delta_k(f)
\end{equation}
where $\bfM_{\gamma,F}^k$ are differential operators defined in 
\eqref{def of Mk gamma F}
and where the remainder term $R_k^\Delta(f)$ 
is defined in \eqref{def of RkDelta}.
Moreover, the differential operators $\bfM_{\gamma,F}^k$
are of order $\leq k$ in each of the variables $h_1,\ldots,h_d$.
Also, the remainder $R_k^\Delta(f)$ satisfies
the estimate of Proposition \ref{polyest}.
\end{Theorem}

In the case where $f$ is a polynomial, this formula gives rise to
an {\em exact} Euler Maclaurin formula.

\begin{Corollary}\labell{polys}
Let $p$ be a polynomial on $\R^n,$ and choose $k \geq  {\rm deg~} p+ n + 1.$
Then 
$${\sum_{\Delta \cap \Z^n}}' p = 
\left.\sum_F \sum_{\gamma \in \Gamma_F^\flat}
       \bfM_{\gamma,F}^k \int_{\Delta(h)} p(x) dx \right|_{h=0}
.$$
\end{Corollary}

\begin{proof}[Proof of Theorem \ref{syms}:]
Let $\chi$ be a smooth function on $\R^n$ which is equal to one
on some open ball about the origin that contains $\Delta$
and which is supported in some larger ball, say, of radius $R$.
Define $\chi_\lambda(x) := \chi(x/\lambda)$ for all $\lambda \geq 1.$
Then the function 
$$ f_\lambda:= f \chi_\lambda$$ 
is a smooth compactly supported function on $\R^n$.
We apply Theorem \ref{main} to obtain
\begin{equation} \labell{EM f lambda}
{\sum_{\Delta \cap \Z^n}}' f_\lambda = 
\left.\sum_F \sum_{\gamma \in \Gamma_F^\flat}
       \bfM_{\gamma,F}^k \int_{\Delta(h)} f_\lambda(x) dx \right|_{h=0}
       + R^\Delta_k(f_\lambda),
\end{equation}
which is valid for any $k \geq 1$.  

Since $f_\lambda$ equals $f$ on a neighborhood of $\Delta$
if $\lambda \geq 1$, the left hand side and the first of the 
two summands on the right hand side in \eqref{EM f lambda}
are equal to the corresponding terms in \eqref{EM formula for symbols}.
Thus, to deduce \eqref{EM formula for symbols} from \eqref{EM f lambda},
it suffices to prove the following claim: 
\begin{equation} \labell{limit claim}
\lim_{\lambda\to \infty}R^\Delta_k(f_\lambda) = R^\Delta_k(f)
\quad \text{ if } \quad k \geq N + n + 1 .
\end{equation}

To prove this claim, we apply the estimate in Proposition \ref{polyest}
to the difference $R_k^\Delta(f) - R_k^\Delta(f_\lambda) 
= R_k^\Delta(f(1-\chi_\lambda))$.  We expand each of the derivatives 
appearing in the estimate into a finite sum of products of derivatives 
of $f$ and derivatives of $1-\chi_\lambda$.  
The leading term has the form $g (1 - \chi_\lambda)$
where $g$ is a derivative of $f$ of order $q$ with $q \geq k$.  
Because $f$ is a symbol of order $N$, the function $g$ is dominated 
by a constant times the function $x \mapsto (1+|x|)^{N-q}$.
Because $q \geq k \geq N+n+1$, the function $g$ is in $L_1$.
It follows that the $L_1$ norm of $g (1-\chi_\lambda)$
converges to zero as $\lambda \to \infty$.
Each of the remaining terms in the expansion has the form
\begin{equation} \labell{the term}
 \lambda^{-s} g(\cdot) \tilde{\chi} (\cdot/\lambda) 
\end{equation}
where $g$ is a derivative of $f$ of order $q$
and $\tilde{\chi}$ is a derivative of $1-\chi$ of order $s$,
and where $s \geq 1$ and $s+q \geq k$.
The function $g$ is dominated by a constant times the 
function $x \mapsto (1+|x|)^{N-q}$; the function $\tilde{\chi}$ is bounded;
the function $\tilde{\chi} (\cdot/\lambda)$ is supported on the ball 
of radius $\lambda R$ about the origin.
Hence, the $L_1$ norm of the product $g(\cdot) \tilde{\chi} (\cdot/\lambda)$ 
is bounded by a constant times $\lambda^{N-q+n}$ if $N-q+n$ is non-negative 
and by a constant otherwise.  
Since $s \geq 1$ and $s+q \geq k \geq N+n+1$, the $L_1$ norm
of the term \eqref{the term} 
is bounded by some constant multiple of $\lambda^{-1}$.  
Letting $\lambda \to \infty$, we see that each of the terms in the estimate
approaches zero as $\lambda \to \infty$.
This implies \eqref{limit claim}.
Theorem \ref{syms} follows.
\end{proof}


\begin{thebibliography}{GDR}

\bibitem[BP]{BP}
A.\ Barvinok and J.\ E.\ Pommersheim,
\emph{An algorithmic theory of lattice points in polyhedra},
New perspectives in algebraic combinatorics (Berkeley, CA, 1996--97),
91--147, Math.\ Sci.\ Res.\ Inst.\ Publ., \textbf{38},
Cambridge Univ.\ Press, Cambridge, 1999.

\bibitem[BDR]{BDR}
M.~Beck, R.~Diaz, and S.~Robins,
\emph{The Frobenius problem, rational polytopes, and Fourier-Dedekind sums},
J.~Number Theory \textbf{96} (2002), no.~1, 1--21.

\bibitem[Bo]{Bo}
N.\ Bourbaki, \emph{Functions d'Une Variable R\'{e}ele},
Chapitre VI: D\'{e}veloppments Tayloriens G\'{e}neralis\'{e}s,
Formule Summatoire d'Euler-Maclaurin (1951).

\bibitem[BV]{BV}
M.~Brion and M.~Vergne, \emph{Lattice points in simple polytopes},
Jour.\ Amer.\ Math.\ Soc.\ \textbf{10} (1997), 371--392.

\bibitem[CS1]{CS:bulletin}
S.~E.~Cappell and J.~L.~Shaneson, \emph{Genera of algebraic varieties
and counting lattice points}, Bull.\ A.~M.~S.\ \textbf{30} (1994), 62--69.

\bibitem[CS2]{CS:EM}
S.~E.~Cappell and J.~L.~Shaneson, \emph{Euler-Maclaurin expansions
for lattices above dimension one},
C.\ R.\ Acad.\ Sci.\ Paris Sér.\ I Math.\ \textbf{321} (1995), 885--890.

\bibitem[CS3]{CS:private}
S.~E.~Cappell and J.~L.~Shaneson, unpublished.

\bibitem[Da]{Da}
V.~I.~Danilov, \emph{The geometry of toric varieties},
Russ.\ Math.\ Surv.\ \textbf{33} (1978) no.~2, 97--154.

\bibitem[DR]{DR} R.\ Diaz and S.\ Robins, \emph{The Ehrhart Polynomial
of a Lattice Polytope}, Ann.\ Math.\ (2) \textbf{145} (1997), 
no.~3, 503--518,
and \emph{Erratum: ``The Ehrhart polynomial of a lattice polytope"},
Ann.\ Math.\ (2) \textbf{146} (1997), no.~1, 237.

\bibitem[Eu]{Eu}
L.~Euler, Commentarii Acad.\ Petrop.\ \textbf{6} (1732--3)
and \textbf{8} (1736).

% \bibitem[GV]{GV}
% A.~N.~Varchenko and I.~M.~Gelfand, \emph{Heaviside functions of
% a configuration of hyperplanes.} (Russian), Funktsional.\ Anal.\
% i Prilozhen.\ \textbf{21} (1987), no.~4, 1--18, 96.
% English translation: Functional Anal.\ Appl.\ \textbf{21} (1987),
% no.~4, 255--270.

\bibitem[Gi]{Gibbs} J. W. Gibbs. ``Fourier Series.'' Nature 59, 200 and 606, 
1899.

\bibitem[Gu1]{Gu:book}
V.~Guillemin, \emph{Moment maps and combinatorial invariants
of Hamiltonian $T^n$-spaces}, Progress in Math.\ \textbf{122}, 
Birkhauser, 1994.

\bibitem[Gu2]{Gu}
V.~Guillemin, \emph{Riemann-Roch for toric orbifolds}, J.\ Diff.\ Geom.\
\textbf{45} (1997), 53--73.

\bibitem[GS]{GS}
V.~Guillemin and S.~Sternberg, \emph{Geometric quantization
and multiplicities of group representations}, Invent.\ Math.\ \textbf{67}
(1982), 515--538.

\bibitem[Har]{H}
G.\ H.\ Hardy, \emph{Divergent Series},
Second Edition, 1991 (unaltered), Chelsea Pub.\ Col, New-York, N.Y.

\bibitem[Hat]{Ha}
A.~Hattori and M.~Masuda, \emph{Theory of multi-fans},
Osaka J.~Math.\ \textbf{40} (2003), 1--68.

\bibitem[Ho]{horm} L. Hormander.  Lectures on nonlinear hyperbolic
differential equations.  New York:  Springer, 1997.

\bibitem[KK]{KK}
J.~M.~Kantor and A.~G.~Khovanskii,
\emph{Une application du th\'eor\`eme de Riemann-Roch combinatoire
au polyn$\hat{o}$me d'Ehrhart des polytopes entiers de $R\sp d$},
C.\ R.\ Acad.\ Sci\. Paris S\'er.\ I Math.\ \textbf{317} (1993),
no.\ 5, 501--507.

\bibitem[KSW]{PNAS}
Y. Karshon, S. Sternberg, and J. Weitsman.
\emph{The Euler-Maclaurin formula for simple integral polytopes},
Proc.\ Nat.\ Acad.\ Sci.\ \textbf{100} no.~2 (2003), 426--433.  

\bibitem[Kh1]{Kh1} A.~G.~Khovanski,
\emph{Newton polyhedra and toroidal varieties},
Funktsional.\ Anal.\ i Prilozhen.\ \textbf{11} (1977), no.~4, 56--64.
English tranlation:
Func.\ Anal.\ Appl.\ \textbf{11} (1977), no.~4, 289--296 (1978).

\bibitem[Kh2]{Kh2} A.~G.~Khovanski,
\emph{Newton polyhedra and the genus of complete intersections},
Funktsional.\ Anal.\ i Prilozhen.\ \textbf{12} (1978), no.~1, 51--61.
English tranlation:
Func.\ Anal.\ Appl.\ \textbf{12} (1978), no.~1, 38--46 (1978).

\bibitem[KP1]{KP1} A.~G.~Khovanskii and A.~V.~Pukhlikov,
\emph{Finitely additive measures of virtual polytopes},
Algebra and Analysis \textbf{4} (1992), 161--185;
translation in St.\ Petersburg Math.\ J.\ \textbf{4} (1993), 
no.~2, 337--356.

\bibitem[KP2]{KP2} A.~G.~Khovanskii and A.~V.~Pukhlikov,
\emph{The Riemann-Roch theorem for integrals and sums of quasipolynomials
on virtual polytopes}, Algebra and Analysis \textbf{4} (1992), 188--216,
translation in St.\ Petersburg Math.\ J.\ (1993), no.~4, 789--812.

\bibitem[Kn]{Kn} K.\ Knopp, \emph{Theory and application of infinite series},
Translated from the
Second German Edition by Miss R.~C.~Young L.\`{e}s Sc., Blackie and Son
Limited, London and Glasgow (1928), Chapter XIV.
 
\bibitem[L]{lawrence}
J.\ Lawrence, \emph{Polytope volume computation},
Math.\ Comp.\ \textbf{57} (1991), no.\ 195, 259--271.

\bibitem[Ma]{Ma}
C.~Maclaurin, \emph{A treatise of Fluxions}, Edinburgh (1742).

\bibitem[Md]{Md} L.~J.~Mordell, \emph{Lattice points in a tetrahedron
and generalized Dedekind sums}, J.~Indian Math.\ Soc.\ (N.S.) 
\textbf{15} (1951), 41--46.

\bibitem[Mo]{Mo}
R.\ Morelli, \emph{Pick's theorem and the Todd class of a toric variety},
Adv.\ Math.\ \textbf{100} (1993), no.\ 2, 183--231.

\bibitem[Od]{Oda}
T.\ Oda, \emph{Convex bodies and algebraic geometry}, Ergebnisse der
Mathematik, Springer, 1988.

\bibitem[Pi]{pick} G. Pick, \emph{Geometrisches zur Zahlentheorie},
Sitzenber.~Lotos (Prague) \textbf{19} (1899), 311--319.

\bibitem[Poi]{Poi}
S.~D.~Poisson, M\'{e}moires Acad.\ sciences Inst.\ France \textbf{6} (1823).

\bibitem[Pom]{Po}
J.\ E.\ Pommersheim,
\emph{Toric varieties, lattice points and Dedekind sums},
Math.\ Ann.\ \textbf{295} (1993), 1--24.

\bibitem[S]{S} J. Shaneson, \emph{Characteristic Classes, Lattice
 Points, and Euler-MacLaurin Formulae}, Proceedings of the International
 Congress of Mathematicians, Zurich, 1994.  Basel:  Birkhauser Verlag,
1995.

\bibitem[T]{T} R.\ M.\ Trigub,
\emph{A Generalization of the Euler-Maclaurin Formula},
Mathematical Notes \textbf{61}, no.~2, 1997, p.~253--257,
translated from Matematicheskie Zametki \textbf{61}, no.~2,
p.~312--316.

\bibitem[V]{V} A.~N.~Varchenko,
\emph{Combinatorics and topology of the arrangement of affine hyperplanes
in the real space.}  (Russian)
Funktsional.\ Anal.\ i Prilozhen.\ \textbf{21} (1987), no.~1, 11--22.
English translation: Functional Anal.\ Appl.\ \textbf{21} (1987), 
no.~1, 9--19.

\bibitem[W]{Wi} W. Wirtinger \emph{Einige Anwendungen der
Euler-Maclaurin'schen Summenformel, insbesondere auf eine Aufgabe
von Abel} Acta. Math. \textbf{26} 
(1902)  255- 271.


\end{thebibliography}
\end{document}